\documentclass[12pt]{amsart}
\usepackage{amsxtra,latexsym,amssymb, amsmath,mathrsfs}
\usepackage{epsfig,rotate,amsthm}
\usepackage{hyperref}

\def\qed{{\hfill $\Box$}}

\def\N{{\mathbb N}}
\def\Z{{\mathbb Z}}
\def\C{{\mathbb C}}

\def\K{{\mathbb K}}

\def\A{ {\mathcal{A}_{n}^{\overline{q},\Lambda}}(\K)}
\def\B{A_{n}^{\overline{q}, \Gamma}(\K)}
\def\C{\mathcal{A}_{m}^{\overline{q}^{\prime}, \Lambda^{\prime}}(\K)}
\def\D{A_{m}^{\overline{q}^{\prime},\Gamma}(\K)}
\def\E{ {\mathcal{A}_{n}^{\overline{q},\Lambda}}(\K[t])}
\def\F{A_{n}^{\overline{q}, \Gamma}(\K[t])}
\def\G{\mathcal{A}_{m}^{\overline{q}^{\prime}, \Lambda^{\prime}}(\K[t])}
\def \H{A_{m}^{\overline{q}^{\prime}, \Gamma^{\prime}}(\K[t])}
\theoremstyle{theorem}
\newtheorem{thm}{Theorem}[section]
\newtheorem{cor}{Corollary}[section]
\newtheorem{prop}{Proposition}[section]

\newtheorem{lem}{Lemma}[section]
\theoremstyle{definition}
\newtheorem{defn}{Definition}[section]
\theoremstyle{remark}

\newtheorem{rem}{\bf Remark}[section]

\begin{document}
\title[``Symmetric" Multiparameter Quantized Weyl Algebras]{Automorphisms for Some ``symmetric" Multiparameter Quantized Weyl Algebras and Their Localizations}
\dedicatory{To Professor Yingbo Zhang on the occasion of her 70th birthday}
\author[X. Tang]{ Xin Tang}
\address{
Department of Mathematics \& Computer Science\\
Fayetteville State University\\
1200 Murchison Road, Fayetteville, NC 28301}
\email{xtang@uncfsu.edu} 

\thanks{This research is partially supported by a research mini-grant funded by the HBCU STEM Master's Program at Fayetteville State University}
\date{\today}
\begin{abstract}
In this paper, we study the algebra automorphisms and isomorphisms for a family of ``symmetric" multiparameter quantized Weyl algebras $\A$ and some related algebras in the generic case. First, we compute the Nakayama automorphism for $\A$ and give a necessary and sufficient condition for $\A$ to be Calabi-Yau. We also prove that $\A$ is cancellative. Then we determine the automorphism group for $\A$ and its polynomial extension $\E$. As an application, we solve the isomorphism problem for $\{\A\}$ and $\{\E\}$. Similar results will be established for the Maltisiniotis multiparameter quantized Weyl algebra $\B$ and its polynomial extension $\F$. In addition, we prove a quantum analogue of the Dixmier conjecture for a simple localization $(\A)_{\mathcal{Z}}$ of $\A$. Moreover, we will completely determine the algebra automorphism group for $(\A)_{\mathcal{Z}}$.
\end{abstract}

\maketitle 
\section*{Introduction}
Ever since it was introduced by Maltisiniotis in \cite{Maltsiniotis}, the multiparameter quantized Weyl algebra $\B$ has been extensively studied in the literature \cite{AD, GZ, Goodearl, FKK, Jordan, Rigal}. When $n=1$, $\B$ is the rank-one quantized Weyl algebra $A_{1}^{q}(\K)$. Recall that $A_{1}^{q}(\K)$ is a $\K-$algebra generated by $x, y$ subject to the relation $xy-qyx=1$, whose prime ideals of $A_{1}^{q}(\K)$ were classified in \cite{Goodearl}. The automorphism group of $A_{1}^{q}(\K)$ was completely determined in \cite{AD}. There has been further research on the Maltisiniotis multiparameter quantized Weyl algebras $\B$ of higher ranks. For instance, Rigal determined the prime spectrum and the automorphism group for $\B$ in the generic case \cite{Rigal}. Later on, Rigal's result has been improved by Goodearl and Hartwig in \cite{GH}, using a result due to Jordan \cite{Jordan} that $\B$ has a simple localization when none of $q_{i}$ is a root of unity. As an application, they have further solved the isomorphism problem for the family of Maltsiniotis multiparameter quantized Weyl algebras $\{\B\}$. 

Another family of multiparameter quantized Weyl algebras $\A$ with symmetric relations has been investigated in the literature \cite{AJ}. We will refer to this family of algebras as the multiparameter quantized Weyl algebras of ``symmetric type" or ``symmetric" multiparameter quantized Weyl algebras. When $n=1$, the multiparameter quantized Weyl algebra $\A$ is also isomorphic to the rank-one quantized Weyl algebra $A_{1}^{q}(\K)$. Note that $\A$ and $\B$ are closely related in many aspects. In particular, it was proved in \cite{Jordan, AJ} that $\A$ and $\B$ have isomorphic simple localizations if none of the major parameters $q_{1}, \cdots, q_{n}$ is a root of unity. The prime ideals of $\A$ were classified in \cite{AJ} in the generic case. When $\lambda_{ij}=1$, the algebra $\A$ is isomorphic to the tensor product $A_{1}^{q_{1}}(\K)_{\K}\otimes \cdots \otimes_{\K} A_{1}^{q_{n}}(\K)$, whose automorphism group has recently been settled in the works \cite{CPWZ1, CPWZ2} for $q_{i}\neq 1$, using the methods of discriminants and Mod p. 

Since $\A$ can be presented as an iterated skew polynomial algebra, $\A$ is a twisted (or skew) Calabi-Yau algebra by the result in \cite{LWW}. Thus, it is of interest to further determine when $\A$ is indeed Calabi-Yau \cite{Ginzburg}. We will compute the Nakayama automorphism for $\A$ using the methods as developed in \cite{LWW, GY} and establish a necessary and sufficient condition for $\A$ to be a Calabi-Yau algebra. We will also prove that $\A$ is universally cancellative when none of $q_{i}$ is a root of unity in the sense of \cite{BZ}. Similar results will be established for the Maltisiniotis multiparameter quantized Weyl algebra $\B$.

The algebra $\A$ also fits into the class of generalized Wey algebras \cite{Bavula}. The isomorphism problem for some rank-one generalized Weyl algebras has been studied in \cite{BJ, RS}. In general, it would be interesting to determine the automorphism group for $\A$ and obtain an isomorphism classification for the family $\{\A\}$. One of the main objectives of this paper is to completely determine the algebra automorphism group for $\A$ and solve the isomorphism problem in the case where none of $q_{1}, \cdots, q_{n}$ is a root of unity. We will use the fact that $\A$ has a simple localization $(\A)_{\mathcal{Z}}$ with respect to the submonoid $\mathcal{Z}$ of $\A$ generated by some normal elements. We will also study the automorphisms and solve the isomorphism problem for the polynomial extensions $\{\E\}$ and $\{\F\}$. Since we don't put any conditions on the parameters $\lambda_{ij}$ and $\gamma_{ij}$, we are not able to employ the method used in \cite{CPWZ2}, due to the lack of information on the discriminants in the root of unity case, despite the recent progress in \cite{CKZ, GKM, LY, NTY}. We will follow the approach used \cite{Rigal} and \cite{GH} in the non-root of unity case. 

When $q=1$ and $n=1$, both $\A$ and $\B$ are indeed isomorphic to the classical first Weyl algebra $A_{1}(\K)$, which is generated by $x,y$ subject to the relation $xy-yx=1$. The automorphism group of $A_{1}(\K)$ was first determined by Dixmier in \cite{Dixmier}, where Dixmier also asked the question whether each $\K-$algebra endomorphism of the $n-$th Weyl algebra $A_{n}(\K)$ is an algebra automorphism when the base field $\K$ is of zero characteristic. Later on, Dixmier's question has been referred to as the Dixmier conjecture. The Dixmier conjecture has been proved to be stably equivalent to the Jacobian conjecture \cite{Keller, BK,Tsuchimoto}. There have been several works studying a quantum analogue of the Dixmier conjecture for some quantum algebras \cite{Backelin, Richard, KL, KL1, Tang}. In particular, it has recently been proved in \cite{KL1} that each $\K-$algebra endomorphism of a simple localization of $A_{1}^{q}(\K)_{\K}\otimes \cdots \otimes_{\K} A_{1}^{q}(\K)$ is indeed an algebra automorphism when $q$ is not a root of unity. It would be of great interest to establish such a quantum analogue for the simple localization $(\A)_{\mathcal{Z}}$ in the generic case. Another main objective of this paper is to address such a problem. We will prove that each $\K-$algebra endomorphism of $(\A)_{\mathcal{Z}}$ is indeed an algebra automorphism under the condition that the parameters $q_{i}$ are independent. We will be able to determine the $\K-$algebra automorphism group for $(\A)_{\mathcal{Z}}$ in this case.

The paper is organized as follows. In Section $1$, we recall the definition of $\A$ and establish some of its basic properties. Then we determine the height-one prime ideals, the normal elements and the center of $\A$. As an application, we prove that $\A$ is universally cancellative. We also compute the Nakayama automorphism of $\A$. Similar results will be established for $\B$ as well. In Section $2$, we completely determine the automorphism group for $\A$ and solve the isomorphism problem. We also study the automorphism group and the isomorphism problem for $\E$ and $\F$. In Section $3$, we prove a quantum analogue of the Dixmier conjecture for $(\A)_{\mathcal{Z}}$ and describe its automorphism group.

\section{Multiparameter Quantized Weyl Algebras of ``Symmetric Type"}
In this section, we first recall the definitions of the Maltisiniotis multiparameter quantized Weyl algebra $\B$ and the ``symmetric" multiparameter quantized Weyl algebra $\A$. Then we will establish some basic properties for $\A$. In particular, we will determine all the height-one prime ideals of $\A$, and describe the normal elements and the center for $\A$ in the case where none of $q_{i}$ is a root of unity. As an application, we show that $\A$ is universally cancellative. We will also verify that $\A$ is a twisted Calabi-Yau algebra, and compute its Nakayama automorphism, and determine when $\A$ is a Calabi-Yau algebra. Similar results will be established for the Maltsiniotis multiparameter quantized Weyl algebra $\B$ as well.

Let $\K$ denote a field and $n$ be any positive integer and set $\overline{q}=(q_{1}, \cdots, q_{n})$ with $q_{i}\in \K^{\ast}$ for $i=1, \cdots, n$. Let $\Gamma=(\gamma_{ij})$ be a multiplicatively skew-symmetric $n\times n$ matrix with $\gamma_{ij}\in \K^{\ast}$. Recall that the Maltisiniotis multiparameter quantized Weyl algebra $\B$ is defined to be a $\K-$algebra generated by $x_{1}, y_{1}, \cdots, x_{n}, y_{n}$ subject to the following relations:
\begin{eqnarray*}
y_{i}y_{j}&=&\gamma_{ij} y_{j}y_{i},\quad \forall i, j \in \{1, 2, \cdots, n\};\\
x_{i}x_{j}&=&q_{i}\gamma_{ij}x_{j}x_{i},\quad \forall 1\leq i<j\leq n;\\
x_{i}y_{j}&=&\gamma_{ji}y_{j}x_{i}, \quad \forall 1\leq i<j\leq n;\\
x_{i}y_{j}&=&q_{i}\gamma_{ji}y_{j}x_{i}, \quad\forall 1\leq j<i\leq n;\\
x_{i}y_{i}-q_{i}y_{i}x_{i}&=& 1+\sum_{k=1}^{i-1}(q_{k}-1)y_{k}x_{k}, i=1, \cdots, n.
\end{eqnarray*}

Now we recall the definition of the so-called alternative (or ``symmetric") multiparameter quantized Weyl algebra $\A$ from \cite{AJ}.
\begin{defn}
Let $\overline{q}=(q_{1}, \cdots, q_{n})$ with $q_{i}\in \K^{\ast}$ and $\Lambda=(\lambda_{ij})_{1\leq i, j\leq n}$ be a square matrix in $M_{n}(\mathbb{K}^{\ast})$ with $\lambda_{ij}=\lambda_{ji}^{-1}$ and $\lambda_{ii}=1$. The multiparameter quantized Weyl algebra $\A$ of ``symmetric type" is defined to be a $\K-$algebra generated by $x_{1}, y_{1}, \cdots, x_{n}, y_{n}$ subject to the relations:
\begin{eqnarray*}
x_{i}x_{j}=\lambda_{ij}x_{j}x_{i},\quad y_{j}y_{i}=\lambda_{ji}y_{i}y_{j}, \quad \forall 1\leq i<j\leq n;\\
x_{i}y_{j}=\lambda_{ji}y_{j}x_{i},\quad x_{j}y_{i}=\lambda_{ij}y_{i}x_{j}, \quad \forall 1\leq i<j\leq n;\\
x_{i}y_{i}-q_{i}y_{i}x_{i}=1, \quad\quad\quad \forall 1\leq i\leq n.
\end{eqnarray*}
\end{defn}
\qed

As we can see, the algebra $\A$ has more symmetric relations compared with the Maltisiniotis multiparameter quantized Weyl algebra $\B$. Next we recall a few well-known results for $\A$. First of all, we have the following proposition.
\begin{prop}
The algebra $\A$ is a Noetherian domain with a Gelfand-Kirillov dimension of $2n$. Indeed, it has a $\K-$basis given as follows:
\[
\mathcal{B}=\{y_{1}^{a_{1}}x_{1}^{b_{1}}\cdots y_{n}^{a_{n}}x_{n}^{b_{n}}| a_{1}, b_{1}, \cdots, a_{n}, b_{n} \in \mathbb{Z}_{\geq 0}\}.
\]
\end{prop}
{\bf Proof:} From the definition of $\A$, it is easy to see that $\A$ can be presented as an iterated skew polynomial algebra as follows:
\[
\A=\K[y_{1}][x_{1}; \tau_{2}, \delta_{2}][y_{2};\tau_{3}, \tau_{4}, \delta_{4}]\cdots [y_{n};\tau_{2n-1},\delta_{2n-1}][x_{n}; \tau_{2n}, \delta_{2n}].
\]
As a result, $\A$ is a (both left and right) Noetherian domain and $\A$ has a Gelfand-Kirillov dimension of $2n$. It also has the $\K-$basis as stated.
\qed

Let $D$ be any ring and $\sigma=(\sigma_{1}, \cdots, \sigma_{n})$ be a set of commuting automorphisms of $D$ and $a=(a_{1}, \cdots, a_{n})$ be a set of (non-zero) elements in the center of $D$ such that $\sigma_{i} (a_{j})=a_{j}$ for $i\neq j$. The study of the {\it generalized Weyl algebra} $A=D(\sigma, a)$ of degree $n$ over the base ring $D$ was initiated by Bavula. For the detailed definition of generalized Weyl algebras and their properties, we refer the readers to \cite{Bavula} and the references therein. Next, we show that $\A$ can be realized as a generalized Weyl algebra. 

For $i=1, \cdots, n$, let us set $z_{i}=x_{i}y_{i}-y_{i}x_{i}$. It is easy to verify that
\[
z_{i}=1+(q_{i}-1)y_{i}x_{i}=q_{i}^{-1}(1+(q_{i}-1)x_{i}y_{i})
\]
and
\begin{eqnarray*}
z_{j}y_{i}=y_{i}z_{j}\,\, \text{if}\, i\neq j,\quad\quad z_{j}y_{j}=q_{j}y_{j}z_{j};\\
z_{j}x_{i}=x_{i}z_{j}\, \,\text{if}\, i\neq j, \quad\quad z_{j}x_{j}=q_{j}^{-1}x_{j}z_{j};\\
z_{i}z_{j}=z_{j}z_{i}.
\end{eqnarray*}
In particular, each $z_{i}$ is a normal element of $\A$ in the sense that $z_{i}\A=\A z_{i}$ for $i=1, \cdots, n$.

\begin{prop}
The algebra $\A$ can be presented as a generalized Weyl algebra of degree $n$. As a result, $\A$ has a $\K-$basis $\mathcal{B}$ consisting of elements which are monomials in $x_{i}, y_{i}, z_{i}$ given as follows:
\[
\mathcal{B}=\{z_{1}^{a_{1}}\cdots z_{n}^{a_{n}}x_{1}^{b_{1}}\cdots x_{n}^{b_{n}}, z_{1}^{a_{1}}\cdots z_{n}^{a_{n}} y_{1}^{c_{1}}\cdots y_{n}^{c_{n}}\mid a_{i}, b_{i}, c_{i}\in \Z_{\geq 0}\}.
\] 
\end{prop}
{\bf Proof:} We can set the base ring as $D=\K[z_{1}, \cdots, z_{n}]$, which is the polynomial $\K-$algebra in $z_{1}, \cdots, z_{n}$. In addition, for $i=1,\cdots, n$, we can choose $a_{i}=\frac{z_{i}-1}{q_{i}-1}$ and define $\sigma_{i}$ by $\sigma_{i}(z_{j})=q_{i}^{\delta_{ij}} z_{j}$ for $j=1, \cdots, n$. Now it is straightforward to check that $\A=D(\sigma, a)$ is a generalized Weyl algebra of degree $n$ over the base ring $D$. By a basic property of the generalized Weyl algebras \cite{Bavula}, we can see that $\A$ has the desired $\K-$basis. 
\qed

Let $\mathcal{Z}$ denote the submonoid of $\A$ generated by the normal elements $z_{1}, \cdots, z_{n}$. Obviously, $\mathcal{Z}$ is an Ore set. Thus we can localize $\A$ with respect to $\mathcal{Z}$. Let $(\A)_{\mathcal{Z}}$ denote the localization of $\A$ with respect to $\mathcal{Z}$. Similarly, one can set $z_{i}^{\prime}=x_{i}y_{i}-y_{i}x_{i}$ in $\B$ for $i=1, \cdots, n$ and denote by $\mathcal{Z}^{\prime}$ the submonoid of $\B$ generated by $z_{1}^{\prime}, \cdots, z_{n}^{\prime}$. The following result has been established in \cite{AJ}.
\begin{thm}
If none of $q_{i}$ is a root of unity, then the localization $(\mathcal{A}_{n, \mathbb{K}}^{Q, \Lambda})_{\mathcal{Z}}$ is a simple algebra.
\end{thm}
{\bf Proof:} Note that it has been proved in \cite{AJ} that the localization $(\A)_{\mathcal{Z}}$ of $\A$ is isomorphic to the corresponding localization $(\B)_{\mathcal{Z}^{\prime}}$ of the Maltsiniotis quantized Weyl algebra $\B$. By {\bf Theorem 3.2} in \cite{Jordan}, the localization $(\B)_{\mathcal{Z}^{\prime}}$ is a simple algebra when none of $q_{i}$ is a root of unity. Thus, the result follows.
\qed

The prime ideals of the ``symmetric" multiparameter quantized Weyl algebra $\A$ were completely determined by Akhavizadgan and Jordan in \cite{AJ} under certain conditions on the parameters $q_{i}$ and $\lambda_{ij}$ in the case where $\K$ is algebraically closed. However, these conditions are only needed for describing the complete prime spectrum of $\A$. As we can see, one can still obtain a complete classification of the height-one prime ideals for the algebra $\A$ under the weaker condition that none of $q_{i}$ is a root of unity and $\K$ is any field.

\begin{thm}
If none of $q_{i}$ is a root of unity, then the height-one prime ideals of $\A$ are exactly given as follows:
\[
\{(z_{1}), (z_{2}), \cdots, (z_{n})\}.
\]
\end{thm}

{\bf Proof:} First of all, by {\bf Lemma 3.2} in \cite{AJ}, each two-sided ideal $(z_{i})$ of $\A$ generated by the normal element $z_{i}$ is a completely prime ideal. Note that the proof of {\bf Lemma 3.2} in \cite{AJ} does not need the extra assumption that $\K$ is an algebraically closed field. In addition, by {\bf Corollary 4.1.12} in \cite{MR}, the height of the prime ideal $(z_{i})$ is at most one for $i=1, \cdots, n$. Since $\A$ is a domain, the zero ideal $(0)$ is a prime ideal of $\A$. Thus, the height of $(z_{i})$ is indeed $1$ for $i=1, \cdots, n$.

Since none of the parameters $q_{i}$ is a root of unity, the localization $(\A)_{\mathcal{Z}}$ is a simple algebra. Thus each non-zero prime ideal $P$ of $\A$ has to meet with the Ore set $\mathcal{Z}$. Suppose that we have $z_{1}^{a_{1}}\cdots z_{n}^{a_{n}}\in P$. Thus we have the inclusion: $(z_{1}^{a_{1}}\cdots z_{n}^{a_{n}})\subseteq P$. Since each $z_{i}$ is a normal element, we have that $z_{j}\in P$ for some $j$. As a result, we have the inclusion: $(z_{j})\subseteq P$ for some $j$. If $P$ is of height one, then we have that $P=(z_{j})$ for some $j$. Therefore, we have completed the proof.
\qed

\begin{cor}
If none of $q_{i}$ is a root of unity. Then the algebra $\A$ has a trivial center $\K$. And the set of all normal elements of $\A$ is given as follows:
\[
N=\{a z_{1}^{l_{1}}\cdots z_{n}^{l_{n}}\mid l_{1}, \cdots, l_{n}\geq 0, a\in \K\}.
\] 
\end{cor}
{\bf Proof:} First of all, it is easy to check that each element $a z_{1}^{i_{1}}\cdots z_{n}^{i_{n}}$ in the set $N$ as given above is indeed a normal element for the algebra $\A$. Conversely, let $w\neq 0$ be a normal element of $\A$. Then the ideal $(w)$ of $\A$ generated by the normal element $w$ is a non-zero two-sided ideal of $\A$. Thus the localization $(w)_{\mathcal{Z}}$ of the ideal $(w)$ with respect to the Ore set $\mathcal{Z}$ is a non-zero two-sided ideal of the localization $(\A)_{\mathcal{Z}}$. Since the localization $(\A)_{\mathcal{Z}}$ is a simple algebra when none of $q_{i}$ is a root of unity, we have that $(w)_{\mathcal{Z}}=(\A)_{\mathcal{Z}}$. Thus we have that $1\in (w)_{\mathcal{Z}}$, which implies that $w^{\prime} wz_{1}^{l_{1}}\cdots z_{n}^{l_{n}}=1$ for some $w^{\prime}\in \A$. So $w$ is invertible in the localization $(\A)_{\mathcal{Z}}$. Note that the only invertible elements of $(\A)_{\mathcal{Z}}$ are non-zero scalar multiples of the products of integer powers of the elements $z_{i}$. Since $w\in \A$, we further have that $w=az_{1}^{l_{1}}\cdots z_{n}^{l_{n}}$ for some $a\in \K^{\ast}$ and $l_{1}, \cdots, l_{n}\in \Z_{\geq 0}$. 

Let $c$ be in the center of $\A$. Since any central element of $\A$ is also a normal element of $\A$. Thus $c=az_{1}^{l_{1}}\cdots z_{n}^{l_{n}}$ for some $a\in \K^{\ast}$ and $l_{1}, \cdots, l_{n}\in \Z_{\geq 0}$. Since $c$ is in the center, we have 
\[
cy_{i}=y_{i}c
\]
for $i=1, \cdots, n$. Note that 
\[
cy_{i}=(az_{1}^{l_{1}}\cdots z_{n}^{l_{n}})y_{i}=q_{i}^{l_{i}} y_{i} (az_{1}^{l_{1}}\cdots z_{n}^{l_{n}})=q_{i}^{l_{i}} y_{i} c
\]
for $i=1, \cdots, n$. Thus $q_{i}^{l_{i}}=1$ for $i=1, \cdots, n$. Since $q_{i}$ is not a root of unity, we have $l_{i}=0$ for $i=1, \cdots, n$. Thus $c=a\in \K^{\ast}$. So the center of $\A$ is the base field $\K$.
\qed

The study of Zariski cancellation problem for noncommutative algebras has recently been initiated in the work of Bell and Zhang. We first recall a definition on the cancellation property from \cite{BZ}. 

\begin{defn} 
Let A be a $\K-$algebra.
\begin{enumerate}
\item We call $A$ cancellative if $A[t]\cong B[t]$ for some $\K-$algebra B implies that $A\cong B$.\\

\item We call $A$ strongly cancellative if, for any $d \geq 1$, the isomorphism $A[t_{1}, \cdots, t_{d}] \cong B[t_{1}, \cdots , t_{d}$]
for some $\K-$algebra $B$ implies that $A\cong B$.\\

\item We call $A$ universally cancellative if, for any finitely generated commutative
$\K-$algebra and domain $R$ such that $R/I = \K$ for some ideal $I \subset R$ and any $\K-$algebra $B$, $A\otimes R \cong B\otimes R$ implies that $A\cong B$.
\end{enumerate}
\end{defn}

Note that if an algebra $A$ is universally cancellative, then it is strongly cancellative; and if an algebra $A$ is strongly cancellative, then it is cancellative.

\begin{thm}
If none of $q_{i}$ is a root of unity, then the algebra $\A$ is universally cancellative. As a result, $\A$ is strongly cancellative and cancellative.
\end{thm}
{\bf Proof:} If none of $q_{i}$ is a root of unity, then the algebra $\A$ has a trivial center $\K$. Thus $\A$ is universally cancellative by {\bf Proposition 1.3} in \cite{BZ}. As a result, the algebra $\A$ is strongly cancellative and cancellative 
\qed

Similarly, we have the following result concerning the cancellation property for the Maltisiniotis multiparameter quantized Weyl algebra $\B$.
\begin{thm}
If none of $q_{i}$ is a root of unity, the Maltisiniotis multiparameter quantized Weyl algebra $\B$ is universally cancellative; and thus it is strongly cancellative and cancellative.
\end{thm}
{\bf Proof:} It can be verified in a similar fashion that the center of the Maltisiniotis multiparameter quantum Weyl algebra $\B$ is also $\K$ when none of the parameters $q_{i}$ is a root of unity. Using {\bf Proposition 1.3} in \cite{BZ}, we know that the algebra $\B$ is universally cancellative. As a result, $\B$ is strongly cancellative and cancellative.
\qed

Recall from \cite{LWW, RRZ} that a $\mathbb{K}-$algebra $A$ is called $\mu-$twisted (or skew) {\it Calabi-Yau} of dimension $d$, where $\mu$ is a $\K-$algebra automorphism of $A$ and $d\in \mathbb{Z}_{\geq 0}$, provided that
\begin{enumerate}
\item $A$ is homologically smooth in the sense that as a module over $A^{e}\colon=A\otimes_{\mathbb{K}}A^{op}$, $A$ has a finitely generated projective resolution of finite length.
\item ${\rm Ext}^{i}_{A^{e}}(A, A^{e})\cong 0 \, \text{or}\, A^{\mu}\, \text{if}\, i=d$. 
\end{enumerate}
When the above two conditions hold, the automorphism $\mu$ is called the Nakayama automorphism of $A$. In general, the Nakayama automorphism $\mu$ of an algebra $A$ is unique up to an inner automorphism of $A$. An algebra $A$ is called {\it Calabi-Yau} (or CY, for short) \cite{Ginzburg} if $A$ is twisted Calabi-Yau and the Nakayama automorphism $\mu_{A}$ of $A$ is inner .

It has been proved that any Ore extension (or skew polynomial algebra) $E=A[x;\sigma, \delta]$ is a twisted Calabi-Yau algebra of dimension $d+1$ if $A$ is a twisted Calabi-Yau algebra of dimension $d$ in \cite{LWW}, where a formula for computing the Nakayama automorphism $\mu$ of $E$ is also established. It follows immediately that any iterated Ore extension (or iterated skew polynomial algebra) is a twisted Calabi-Yau algebra. The Nakayama automorphisms of iterated Ore extensions (or iterated skew polynomial algebras) have been further studied in \cite{GY}, where a formula of computing the Nakayama automorphism for any reversible, diagonalized iterated Ore extension is developed. In particular, it is proved in \cite{GY} that the Nakayama automorphism of any symmetric CGL extension is an automorphism defined by the conjugation of a normal element. Since $\A$ is an iterated skew polynomial algebra, $\A$ is indeed a twisted Calabi-Yau algebra. We next compute the Nakayama automorphism for $\A$ and provide a necessary-sufficient condition for $\A$ to be a Calabi-Yau algebra. In addition, we will identify a normal element that defines the Nakayama automorphism $\mu$ of $\A$.

\begin{thm}
The ``symmetric" multiparameter quantized Weyl algebra $\A$ is a twisted Calabi-Yau algebra. The Nakayama automorphism $\mu$ of $\A$ is given as follows:
\[
\mu (x_{i})=q_{i}x_{i},\quad \mu (y_{i})=q_{i}^{-1}y_{i}
\]
for $i=1, \cdots, n$. As a result, the algebra $\A$ is a Calabi-Yau algebra if and only if $q_{i}=1$ for $i=1, \cdots, n$. Moreover, the Nakayama automorphism $\mu$ of $\A$ is the automorphism of $\A$ defined by the conjugation of $z=z_{1}z_{2}\cdots z_{n}$. That is, 
\[
\mu(x_{i})=z^{-1}x_{i}z, \quad\quad \mu(y_{i})=z^{-1}y_{i}z
\]
for $i=1, \cdots, n$. 
\end{thm}
{\bf Proof:} We prove this statement using induction on the index $n$. For simplicity, we will sometimes denote the algebra $\A$ by $A$ instead. When $n=1$, we have that $A=\mathcal{A}_{1}^{q_{1}}(\K)$, the rank-one quantized Weyl algebra, which is a $\K-$algebra generated by $x_{1}$ and $y_{1}$ subject to the only relation: $x_{1}y_{1}-q_{1}y_{1}x_{1}=1$. In this case, we can present the algebra $A$ as an iterated skew polynomial algebra in two different ways. First of all, we have that
\[
A=\K[y_{1}][x_{1}; \sigma_{1}, \delta_{1}]
\]
where $\sigma_{1}(y_{1})=q_{1}y_{1}$ and $\delta_{1} (y_{1})=1$. Applying {\bf Theorem 3.3} in \cite{LWW}, we have that $A$ is a twisted Calabi-Yau algebra and its Nakayama automorphism $\mu$ is given as follows:
\[
\mu (y_{1})=q_{1}^{-1}y_{1}, \quad \mu (x_{1})=ax_{1}+b\]
for some invertible element $a\in \K^{\ast}$ and $b$ in $\K[y_{1}]$. In addition, one can also present the algebra $A$ as an iterated skew polynomial algebra as follows:
\[
A=\K[x_{1}][y_{1}; \sigma^{\prime}_{1}, \delta^{\prime}_{1}]
\]
where $\sigma_{1}^{\prime}(x_{1})=q_{1}^{-1}x_{1}$ and $\delta_{1}^{\prime} (x_{1})=-q_{1}^{-1}$. Thus, $A$ is a twisted Calabi-Yau algebra with its Nakayama automorphism $\mu$ defined by: 
\[
\mu (x_{1})=q_{1}x_{1}, \quad \mu (y_{1})=cy_{1}+d
\]
for some invertible element $c\in \K^{\ast}$ and $d$ in $\K[x_{1}]$. Since the Nakayama automorphism of $A$ is unique up to an inner automorphism of $A$ and the set of invertible elements in $A$ is the set $\K^{\ast}$, we have that $b=d=0$ and $a=q_{1}$ and $c=q_{1}^{-1}$.

Suppose the statement is true for the subalgebra $B$ of $\A$ generated by $x_{1}, \cdots, x_{n-1}, y_{1}, \cdots, y_{n-1}$. That is, the algebra $B$ is a twisted Calabi-Yau algebra with its Nakayama automorphism $\mu$ defined as follows:
\[
\mu (x_{i})=q_{i}x_{i}, \quad \mu (y_{i})=q_{i}^{-1}y_{i}
\]
for $i=1, \cdots, n-1$. 

Next we will present $\A$ as an iterated skew polynomial algebra over the base algebra $B$ in two ways. We will denote the Nakayama automorphism of $\A$ by $\mu$. First of all, we can have that $
\A=B[y_{n}; \sigma_{1}, \delta_{1}][x_{n}; \sigma_{2}, \delta_{2}]$
with 
\[
\sigma_{1}(x_{i})=\lambda_{in} x_{i},\quad \sigma_{1}(y_{i})=\lambda_{ni} y_{i}, \quad \delta_{1}(x_{i})=\delta_{1}(y_{i})=0
\]
for $i=1, \cdots, n-1$ and 
\[
\sigma_{2}(y_{n})=q_{n}y_{n}, \quad \delta_{2} (y_{n})=1
\]
and
\[
\sigma_{2}(x_{i})=\lambda_{ni} x_{i},\quad \sigma_{2}(y_{i})=\lambda_{in} y_{i}, \quad \delta_{2}(x_{i})=\delta_{2}(y_{i})=0
\]
for $i=1, \cdots, n-1$. Since $\delta_{1}=0$, by {\bf Theorem 3.2} and {\bf Remark 3.4} in \cite{LWW}, the Nakayama automorphism $\mu$ of $\A$ is given as follow:
\[
\mu (x_{i})=\lambda_{ni}^{-1} \lambda_{in}^{-1} q_{i}x_{i}=q_{i}x_{i}, \quad \mu(y_{i})=\lambda_{in}^{-1}\lambda_{ni}^{-1} q_{i}^{-1} y_{i}=q_{i}^{-1}y_{i}
\]
for $i=1, \cdots, n-1$ and
\[
\mu(y_{n})=a q_{n}^{-1} y_{n}, \quad \mu(x_{n})=c x_{n}+d
\]
for some $a, c\in \K^{\ast}$ and $d\in B[y_{n};\sigma_{1}, \delta_{1}]$. Second of all, we can have that $\A=B[x_{n}; \sigma^{\prime}_{1}, \delta^{\prime}_{1}][y_{n}; \sigma^{\prime}_{2}, \delta^{\prime}_{2}]$
with 
\[
\sigma^{\prime}_{1}(x_{i})=\lambda_{ni} x_{i},\quad \sigma_{1}(y_{i})=\lambda_{in} y_{i}, \quad \delta^{\prime}_{1}(x_{i})=\delta^{\prime}_{1}(y_{i})=0
\]
for $i=1, \cdots, n-1$ and 
\[
\sigma^{\prime}_{2}(x_{n})=q_{n}^{-1}x_{n}, \quad \delta^{\prime}_{2} (x_{n})=-q_{n}^{-1}
\]
and
\[
\sigma_{2}(x_{i})=\lambda_{in} x_{i},\quad \sigma_{2}(y_{i})=\lambda_{ni} y_{i}, \quad \delta_{2}(x_{i})=\delta_{2}(y_{i})=0
\]
for $i=1, \cdots, n-1$. Therefore, the Nakayama automorphism of $\A$ is given as follows:
\[
\mu (x_{i})=\lambda_{in}^{-1} \lambda_{ni}^{-1} q_{i}x_{i}=q_{i}x_{i}, \quad \mu(y_{i})=\lambda_{ni}^{-1}\lambda_{in}^{-1} q_{i}^{-1} y_{i}=q_{i}^{-1}y_{i}
\]
for $i=1, \cdots, n-1$ and
\[
\mu(x_{n})=a^{\prime} q_{n} x_{n}, \quad \mu(y_{n})=c^{\prime} y_{n}+d^{\prime}
\]
for some $a^{\prime}, c^{\prime}\in \K^{\ast}$ and $d^{\prime} \in B[x_{n};\sigma_{1}^{\prime}, \delta_{1}^{\prime}]$. Note that the subalgebra $B[y_{n}; \sigma_{1}, \delta_{1}]$ of $\A$ can also be constructed as an iterated skew polynomial algebra from the base field $\K$ as follows:
\[
B[y_{n}; \sigma_{1}, \delta_{1}]=\K[y_{n}][x_{1};\sigma_{2}^{\prime\prime\prime}, \delta_{2}^{\prime\prime\prime}][y_{1};\sigma_{3}^{\prime\prime\prime}, \delta_{3}^{\prime\prime\prime}] \cdots[y_{n-1};\sigma_{2n-1}^{\prime\prime\prime}, \delta_{2n-1}^{\prime\prime\prime}].
\]
Using {\bf Theorem 3.2} in \cite{LWW} repeatedly, we have $\mu(y_{n})=q_{n}^{-1}y_{n}+e$ for some $e\in B[y_{n}; \sigma_{1}, \delta_{1}]$. Similarly, we can prove that $\mu(x_{n})=q_{n} x_{n}+f$ for some $f\in B[x_{n}; \sigma^{\prime}_{1}, \delta^{\prime}_{1}]$. Comparing the various expressions of the Nakayama automorphism $\mu$ for the algebra $\A$, we have that $d=d^{\prime}=e=f=0$ and $a^{\prime}q_{n}=c=q_{n}$ and $c^{\prime}=aq_{n}^{-1}=q_{n}^{-1}$. Thus the Nakayama automorphism $\mu$ of $\A$ is given as follows:
\[
\mu (x_{i})=q_{i}x_{i},\quad \mu (y_{i})=q_{i}^{-1}y_{i}
\]
for $i=1, \cdots, n$. Therefore, $\A$ is Calabi-Yau if and only if $q_{i}=1$ for $i=1, \cdots, n$.

It is straightforward to check that $\mu$ is indeed defined by the conjugation of $z=z_{1}\cdots z_{n}$. So we have completed the proof.
\qed

We can determine the Nakayama automorphism $\mu$ for the multiparameter Maltsiniotis quantized Weyl algebra $\B$ in a similar fashion. In particular, we have the following result.
\begin{thm}
The Maltsiniotis multiparameter quantized algebra $\B$ is a twisted {\it Calabi-Yau algebra} with a Nakayama automorphism $\mu$ defined as follows:
\[
\mu (x_{i})=q_{i}^{n+1-i}x_{i},\quad \mu (y_{i})=q_{i}^{-n-1+i}y_{i}
\]
for $i=1, \cdots, n$. As a result, the algebra $\B$ is a Calabi-Yau algebra if and only if each $q_{i}$ is an $(n+1-i)$-th root of unity. In particular, the Nakayama automorphism $\mu$ of $\B$ is the automorphism of $\B$ defined by the conjugation of the normal element $z=(x_{1}y_{1}-y_{1}x_{1})\cdots (x_{n}y_{n}-y_{n}x_{n})$ as follows:
\[
\mu(x_{i})=z^{-1} x_{i} z,\quad \quad \mu(y_{i})=z^{-1} y_{i}z
\]
for $i=1, \cdots, n$.
\end{thm}
{\bf Proof:} The proof is similar to the one for {\bf Theorem 1.5}, and we will not repeat the details here.
\qed

\begin{rem}
The Nakayama automorphism $\mu$ of $\A$ (or $\B$) is completely determined by the parameters $q_{1}, \cdots, q_{n}$, and it is independent of the rest parameters $\lambda_{ij}$ (or $\gamma_{ij}$).
\end{rem}
\begin{rem}
One can also present the algebra $\A$ (or $\B$) as a reversible, diagonalized iterated Ore extension and then apply the formula as developed in \cite{GY} to compute the Nakayama automorphism of $\A$ (or $\B$). Instead, we have closely followed the original idea as established in \cite{LWW} for the purpose of providing more details.
\end{rem}
\qed

\section{Automorphisms and Isomorphisms for ``Symmetric" Quantized Weyl Algebras}
In this section, we determine the automorphism groups for $\A$ and $\E$ and classify the algebras $\A$ and $\E$ up to algebra isomorphisms in the case where none of $q_{1}, \cdots, q_{n}$ is a root of unity. Note that we will not put any condition on the rest parameters $\lambda_{ij}$. We have some similar results on the automorphisms and isomorphisms for $\F$.

First of all, we establish a useful lemma, which is similar to the ones used in \cite{Rigal} and \cite{GK}. We will further set $z_{0}=1$.
\begin{lem}
For any $1\leq i\leq n$, let $a,b \in \A\backslash \K$ such that 
\[
ab=a_{0}z_{0}+\sum_{j=1}^{n} \alpha_{j} z_{j}
\]
with $\alpha_{j}\in \K$ and $\alpha_{i}\neq 0$. Then there exist $\lambda, \mu \in \K^{\ast}$ such that 
\[
a=\lambda y_{i},\quad \quad b=\mu x_{i}
\]
or 
\[
a=\lambda x_{i},\quad\quad b=\mu y_{i}.
\]
\end{lem}
{\bf Proof:} For any $1\leq i\leq n$, we can regard the algebra $\A$ as a $\N-$filtered algebra by assigning the degrees to its generators as follows:
\[
\deg (x_{i})=1=\deg (y_{i}), \quad \deg(x_{j})=\deg (y_{j})=0
\]
for all $j\neq i$. Note that the above assignment of degrees is compatible with the defining relations of $\A$, and the degree of any element in $\A$ is defined to be its degree when considered as an element in the corresponding graded algebra. Since $q_{k}\neq 1$ for $k=1, \cdots, n$, we have that $\deg(z_{j})=0$ for $j\neq i$ and $\deg(z_{i})=2$. Since $ab=\sum_{j=0}^{n}\alpha_{j}z_{j}$ with $\alpha_{i}\neq 0$, we have that $\deg(ab)=2$. Due to the degree consideration, we have the following three possibilities: 
\[
\deg (a)=2,\quad \deg (b)=0;
\]
or
\[
\deg(a)=0,\quad \deg (b)=2;
\]
or 
\[
\deg(a)=\deg(b)=1.
\]
In the first two cases, we have that $a=fx_{i}^{2}+gx_{i}y_{i}+hy_{i}^{2}$ or $b=fx_{i}^{2}+gx_{i}y_{i}+hy_{i}^{2}$ for some elements $f, g\in \A$ whose terms do not involve any positive powers of $x_{i}$ and $y_{i}$ at all. As a result, we should have that either $b$ or $a$ is a scalar. This cannot happen. Thus only the third case is possible. That is, we have 
\[
a=fx_{i}+gy_{i},\quad b=f^{\prime}x_{i}+g^{\prime} y_{i}
\]
for some $f, g, f^{\prime}, g^{\prime}\in \A$ whose terms do not involve any positive powers of $x_{i}$ and $y_{i}$. It is easy to see that we have 
\[
ff^{\prime}=0,\quad gg^{\prime}=0.
\]
Suppose that $f=0$, then $g\neq 0$, which implies that $g^{\prime}=0$. Since $b\neq 0$, we have that $f^{\prime}\neq 0$. So we have 
\[
a=gy_{i}, \quad b=f^{\prime} x_{i}.
\]

Note that $ab$ is a $\K-$linear combination of $z_{1}, \cdots, z_{n}$, where $z_{j}=x_{j}y_{j}-y_{j}x_{j}$ for $j=0, \cdots, n$. Thus we have that both $g$ and $f^{\prime}$ are in $\K^{\ast}$. Moreover, we have 
\[
a=\lambda y_{i}, \quad b=\mu x_{i}.
\]
If $f^{\prime}=0$, then we have that $g^{\prime}\neq 0$. Thus we have that $g=0$. So we have 
\[
a=fx_{i},\quad b=g^{\prime} y_{i}.
\]
As a result, we have the following:
\[
a=\lambda x_{i},\quad b=\mu y_{i}
\]
for some $\lambda, \mu \in \K^{\ast}$. So we have completed the proof.
\qed

\begin{lem}
Let $\varphi$ be a $\K-$algebra automorphism for $\A$. Then for each $i=1, \cdots, n$, we have that $\varphi(z_{i})=(z_{j})$ for some $j$. In particular, we have that 
\[
\varphi(x_{i})=\alpha x_{j},\quad \varphi(y_{i})=\beta y_{j}
\] or 
\[
\varphi(x_{i})=\alpha y_{j},\quad \varphi(y_{i})=\beta x_{j}
\]
for some $j$ and $\alpha, \beta \in \K^{\ast}$.
\end{lem}
{\bf Proof:} Since $\varphi$ is a $\K-$algebra automorphism of $\A$, it permutes all the height-one prime ideals of $\A$. As a result, we have 
\[
\varphi((z_{i}))=(z_{j})
\]
for some $j$. Thus $\varphi(z_{i})=fz_{j}$ for some $f\in \A$. Conversely, we have that $\phi^{-1} ((z_{j}))=(z_{i})$. Thus $\varphi^{-1}(z_{j})=gz_{i}$ for some $g\in \A$. As a result, we have the following
\[
z_{i}=\varphi^{-1} (f) g z_{i}
\]
which implies that $\varphi^{-1}(f) g=1$. Since the only invertible elements of $\A$ are in $\K^{\ast}$, we have that $\phi^{-1}(f)\in \K^{\ast}$ and $g\in \K^{\ast}$. Thus we have that $\varphi(z_{i})=\lambda z_{j}$ for some $\lambda \in \K^{\ast}$. Since $x_{i}y_{i}=\frac{q_{i}z_{i}-1}{q_{i}-1}$, we have the following
\[
\varphi(x_{i}) \varphi(y_{i})=\frac{q_{i}\lambda z_{j}-1}{q_{i}-1}=\frac{q_{i}\lambda}{q_{i}-1}z_{j}-\frac{1}{q_{i}-1}.
\]
By {\bf Lemma 2.1}, we have either $\sigma(x_{i})=\alpha x_{j}$ and $\varphi(y_{i})=\beta y_{j}$, or $\varphi(x_{i})=\alpha y_{j}$ and $\varphi(y_{i})=\beta x_{j}$ for some $\alpha, \beta \in \K^{\ast}$ and $j\in \{1, \cdots, n\}$.
\qed

\begin{thm}
Let $\varphi \in {\rm Aut}_{\K}(\A)$. Then there exist a permutation $\sigma$ of the set $\{1, \cdots, n\}$ and $\alpha_{i}\in \mathbb{K}^{\ast}$ for $i=1, \cdots, n$ such that 
\[
\varphi(x_{i})=\alpha_{i}x_{\sigma(i)},\quad \varphi(y_{i})=\alpha_{i}^{-1} y_{\sigma(i)},\quad q_{i}=q_{\sigma(i)},
\]
or 
\[
\varphi(x_{i})=\alpha_{i}y_{\sigma(i)},\quad \quad \varphi(y_{i})=-q_{i}^{-1}\alpha_{i}^{-1} x_{\sigma(i)},\quad q_{i}=q_{\sigma(i)}^{-1}.
\]
\end{thm}

{\bf Proof:} Since $\varphi$ is a $\K-$algebra automorphism of $\A$, we have that $\varphi(x_{i})=\alpha x_{i}$ and $\varphi(y_{i})=\beta y_{i}$ or $\varphi(x_{i})=\alpha y_{j}$ and $\varphi(y_{i})=\beta x_{j}$ for some $j\in \{1, \cdots, n\}$ and $\alpha, \beta \in \K^{\ast}$. Suppose that $\varphi(x_{i})=\alpha x_{j}, \varphi(y_{i})=\beta y_{j}$. Then we have 
\[
\alpha \beta (x_{j}y_{j}-q_{i}y_{j}x_{j})=1, \quad \alpha\beta(x_{j}y_{j}-q_{j}y_{j}x_{j})=\alpha\beta.
\]
As a result, we have that $\alpha\beta (q_{j}-q_{i})y_{j}x_{j}=(1-\alpha\beta)$, which implies that $\alpha\beta =1$, and $q_{i}=q_{j}$. So we have that $\beta=\alpha^{-1}$ and $q_{i}=q_{j}$. 
Suppose that $\varphi(x_{i})=\alpha y_{j}$ and $\varphi(y_{i})=\beta x_{j}$. Then we have the following
\[
\alpha \beta (y_{j}x_{j}-q_{i} x_{j}y_{j})=1,\quad \alpha\beta (x_{j}y_{j}-q_{j}y_{j}x_{j})=\alpha\beta.
\]
As a result, we have that $\alpha\beta (1-q_{i}q_{j})y_{j}x_{j}=(1+q_{i}\alpha\beta)$, which implies that $q_{i}=q_{j}^{-1}$ and $\beta=-q_{i}^{-1} \alpha^{-1}$. In either case, we can define a permutation $\sigma$ of $\{1, \cdots, n\}$ by $\sigma(i)=j$. So the result has been proved.
\qed

The following two corollaries follow directly from {\bf Theorem 2.1} and we state them here without proofs.
\begin{cor}
If $q_{i}\neq q_{j}^{\pm 1}$ for $i\neq j$, then ${\rm Aut}_{\K} (\A)\cong (\K^{\ast})^{n}$.
\end{cor} \qed
\begin{cor}
If $q_{i}=q_{j}$ for all $i, j$, then ${\rm Aut}_{\K}(\A)\cong G \ltimes (\K^{\ast})^{n}$ where $G$ is a finite subgroup of ${\rm Aut}_{\K} (\A)$ consisting of the automorphisms $\A$ defined by the permutations of $\{1, \cdots, n\}$, which are compatible with the defining relations of $\A$.
\end{cor}\qed

Next we are going to give a more explicit description for the automorphism group of $\A$.
\begin{thm}
Let $\varphi \colon \A \longrightarrow \A$ be a $\K-$algebra automorphism of $\A$. Then there exist a partition $\{1, 2, \cdots, n\}=P_{1}\cup P_{2}$ and a permutation $\sigma \colon \{1, 2, \cdots, n\} \longrightarrow \{1, 2, \cdots, n\}$ such that
\begin{enumerate}
\item for any $i, j\in P_{1}$, we have
\[
\lambda_{ij}=\lambda_{\sigma(i) \sigma(j)},\quad\quad q_{i}=q_{\sigma(i)},\quad\quad q_{j}=q_{\sigma(j)}.
\]

\item for $i\in P_{1}, j\in P_{2}$, we have
\begin{eqnarray*}
\lambda_{ij}=\lambda_{\sigma(i)\sigma(j) }^{-1},\quad\quad q_{i}=q_{\sigma(i)}, \quad\quad q_{j}=q_{\sigma(j)}^{-1}.
\end{eqnarray*}
\item for $i, j\in P_{2}$, we have
\begin{eqnarray*}
\lambda_{ij}=\lambda_{\sigma(i)\sigma(j)},\quad\quad
q_{i}=q_{\sigma(i)}^{-1}, \quad\quad q_{j}=q_{\sigma(j)}^{-1}.
\end{eqnarray*}
\end{enumerate}
In particular, the automorphism $\varphi$ is defined as follows:
\begin{eqnarray*}
\varphi(x_{i})=\alpha_{i} x_{\sigma(i)},\quad\quad \varphi(y_{i})=\alpha_{i}^{-1} y_{\sigma(i)},\quad \forall i\in P_{1};\\
\varphi(x_{j})=\alpha_{j} y_{\sigma(j)}, \quad\quad \varphi(y_{j})=-q_{i}^{-1}\alpha_{j}^{-1} x_{\sigma(j)},\quad \forall j\in P_{2}.
\end{eqnarray*}

Conversely, for any given partition $P_{1}\cup P_{2}$ of $\{1, \cdots, n\}$ and any permutation $\sigma$ of $\{1, 2, \cdots, n\}$ satisfying the above conditions, we can define a $\K-$algebra automorphism $\varphi$ of $\A$ as follows:
\begin{eqnarray*}
\varphi(x_{i})=\alpha_{i} x_{\sigma(i)},\quad\quad \varphi(y_{i})=\alpha_{i}^{-1} y_{\sigma(i)},\quad \forall i\in P_{1};\\
\varphi(x_{j})=\alpha_{j} y_{\sigma(j)}, \quad\quad \varphi(y_{j})=-q_{j}^{-1}\alpha_{j}^{-1} x_{\sigma(j)},\quad\forall j\in P_{2}.
\end{eqnarray*} 
\end{thm}

{\bf Proof:} Let $\varphi$ be a $\K-$algebra automorphism of $\A$ and $\sigma$ be the corresponding permutation of $\{1, \cdots, n\}$ associated to $\varphi$. Let us set
\[
P_{1}=\{i\in \{1,\cdots, n\}| \varphi(x_{i})=\alpha_{i} x_{\sigma(i)}, \varphi(y_{i})=\alpha_{i}^{-1} y_{\sigma(i)}\}\]
and
\[
P_{2}=\{i \in \{1, \cdots, n\}| \varphi(x_{i})=\alpha_{i} y_{\sigma(i)}, \varphi(y_{i})=-q_{i}^{-1}\alpha_{i}^{-1} x_{\sigma(i)}\}.\]
Then $P_{1}\cup P_{2}=\{1, \cdots, n\}$ is a partition of the set $\{1, 2, \cdots, n\}$.

Let $i, j\in P_{1}$, then we have the following:
\[
q_{i}=q_{\sigma(i)},\quad q_{j}=q_{\sigma(j)}
\]
and
\[
\varphi(x_{i})=\alpha_{i} x_{\sigma(i)}, \quad \varphi(y_{i})=\alpha_{i}^{-1} y_{\sigma(i)}
\]
and
\[
\varphi(x_{j})=\alpha_{j} x_{\sigma(j)},\quad \varphi(y_{j})=\alpha_{j}^{-1} y_{\sigma(j)}.
\]
From the definition of $\A$, we have the following commuting relations among the generators $x_{i}, x_{j}, y_{i}$, and $y_{j}$:
\begin{eqnarray*}
x_{i}x_{j}=\lambda_{ij} x_{j}x_{i}, \quad y_{i}y_{j}=\lambda_{ij} y_{j}y_{i},\\
x_{i}y_{j}=\lambda_{ji} y_{j}x_{i},\quad y_{i}x_{j}=\lambda_{ji} x_{j}y_{i}.
\end{eqnarray*}
Applying the automorphism $\varphi$ to both sides of the above equations, we will have the following:
\begin{eqnarray*}
x_{\sigma(i)} x_{\sigma(j)}=\lambda_{ij} x_{\sigma(j)} x_{\sigma(i)}, \quad y_{\sigma(i)} y_{\sigma(j)}=\lambda_{ij} y_{\sigma(j)} y_{\sigma(i)},\\
x_{\sigma(i)} y_{\sigma(j)}=\lambda_{ji} y_{\sigma(j)} x_{\sigma(i)},\quad y_{\sigma(i)} x_{\sigma(j)} =\lambda_{ji} x_{\sigma(j)} y_{\sigma(i)}. 
\end{eqnarray*}
By the definition of $\A$, we also have the following:
\begin{eqnarray*}
x_{\sigma(i)} x_{\sigma(j)}=\lambda_{\sigma(i)\sigma(j)} x_{\sigma(j)} x_{\sigma(i)}, \quad y_{\sigma(i)} y_{\sigma(j)}=\lambda_{\sigma(i)\sigma(j)} y_{\sigma(j)} y_{\sigma(i)},\\
x_{\sigma(i)} y_{\sigma(j)}=\lambda_{\sigma(j)\sigma(i)} y_{\sigma(j)} x_{\sigma(i)},\quad y_{\sigma(i)} x_{\sigma(j)} =\lambda_{\sigma(j)\sigma(i)} x_{\sigma(j)} y_{\sigma(i)}. 
\end{eqnarray*}
Therefore, we have $\lambda_{ij}=\lambda_{\sigma(i)\sigma(j)}$ and $q_{i}=q_{\sigma(i)}$ and $q_{j}=q_{\sigma(j)}$ for any $i, j\in P_{1}$. 

Let $i\in P_{1}$ and $j\in P_{2}$, then we have the following:
\[
q_{i}=q_{\sigma(i)},\quad q_{j}=q_{\sigma(j)}^{-1}
\]
and 
\[
\varphi(x_{i})=\alpha_{i} x_{\sigma(i)}, \varphi(y_{i})=\alpha_{i}^{-1} y_{\sigma(i)}
\]
and
\[
\varphi(x_{j})=\alpha_{j} y_{\sigma(j)},\quad \varphi(y_{j})=-q_{i}^{-1} \alpha_{j}^{-1} x_{\sigma(j)}.
\]
Applying the automorphism $\varphi$ to the following equations:
\begin{eqnarray*}
x_{i}x_{j}=\lambda_{ij} x_{j}x_{i}, \quad y_{i}y_{j}=\lambda_{ij} y_{j}y_{i},\\
x_{i}y_{j}=\lambda_{ji} y_{j}x_{i},\quad y_{i}x_{j}=\lambda_{ji} x_{j}y_{i},
\end{eqnarray*}
we have the following equations:
\begin{eqnarray*}
x_{\sigma(i)} y_{\sigma(j)}=\lambda_{ij} y_{\sigma(j)} x_{\sigma(i)}, \quad y_{\sigma(i)} x_{\sigma(j)}=\lambda_{ij} x_{\sigma(j)} y_{\sigma(i)},\\
x_{\sigma(i)} x_{\sigma(j)}=\lambda_{ji} x_{\sigma(j)} x_{\sigma(i)},\quad y_{\sigma(i)} y_{\sigma(j)} =\lambda_{ji} y_{\sigma(j)} y_{\sigma(i)}. 
\end{eqnarray*}

From the definition of $\A$, we also have the following:
\begin{eqnarray*}
x_{\sigma(i)} y_{\sigma(j)}= \lambda_{\sigma(j)\sigma(i)} y_{\sigma(j)} x_{\sigma(i)},\quad y_{\sigma(i)} x_{\sigma(j)}=\lambda_{\sigma(j)\sigma(i)} x_{\sigma(j)}y_{\sigma(i)},\\
x_{\sigma(i)} x_{\sigma(j)}= \lambda_{\sigma(i)\sigma(j)} x_{\sigma(j)} x_{\sigma(i)},\quad y_{\sigma(i)} y_{\sigma(j)}=\lambda_{\sigma(i)\sigma(j)} y_{\sigma(j)} y_{\sigma(i)}.
\end{eqnarray*}
As a result, we have that $\lambda_{ij}=\lambda_{\sigma(j)\sigma(i)}=\lambda_{\sigma(i)\sigma(j)}^{-1}$ for $i\in P_{1}$ and $j\in P_{2}$. 

Let $i, j\in P_{2}$, then we have the following:
\[
q_{i}=q_{\sigma(i)}^{-1}, \quad q_{j}=q_{\sigma(j)}^{-1}
\]
and
\[
\varphi(x_{i})=\alpha_{i} y_{\sigma(i)},\quad \varphi(y_{i})=-q_{i}^{-1}\alpha_{i}^{-1} x_{\sigma(i)}
\]
and
\[
\varphi(x_{j})=\alpha_{j} y_{\sigma(j)},\quad \varphi(y_{j})=-q_{j}^{-1} \alpha_{j}^{-1} x_{\sigma(j)}.
\]

Apply the automorphism $\varphi$ to both sides of the following equations:
\begin{eqnarray*}
x_{i}x_{j}=\lambda_{ij} x_{j}x_{i}, \quad y_{i}y_{j}=\lambda_{ij} y_{j}y_{i},\\
x_{i}y_{j}=\lambda_{ji} y_{j}x_{i},\quad y_{i}x_{j}=\lambda_{ji} x_{j}y_{i},
\end{eqnarray*}
and we have the following:
\begin{eqnarray*}
y_{\sigma(i)} y_{\sigma(j)}=\lambda_{ij} y_{\sigma(j)} y_{\sigma(i)}, \quad x_{\sigma(i)} x_{\sigma(j)}=\lambda_{ij} x_{\sigma(j)} x_{\sigma(i)},\\
y_{\sigma(i)} x_{\sigma(j)}=\lambda_{ji} x_{\sigma(j)} x_{\sigma(i)},\quad x_{\sigma(i)} y_{\sigma(j)} =\lambda_{ji} y_{\sigma(j)} x_{\sigma(i)}. 
\end{eqnarray*}

Additionally, we have the following:
\begin{eqnarray*}
y_{\sigma(i)} y_{\sigma(j)}= \lambda_{\sigma(i)\sigma(j)} y_{\sigma(j)} y_{\sigma(i)},\quad x_{\sigma(i)} x_{\sigma(j)}=\lambda_{\sigma(i)\sigma(j)} x_{\sigma(j)}x_{\sigma(i)},\\
y_{\sigma(i)} x_{\sigma(j)}= \lambda_{\sigma(j)\sigma(i)} x_{\sigma(j)} y_{\sigma(i)},\quad x_{\sigma(i)} y_{\sigma(j)}=\lambda_{\sigma(j)\sigma(i)} y_{\sigma(j)} x_{\sigma(i)}.
\end{eqnarray*}
After comparing these equations, we have that $\lambda_{ij}=\lambda_{\sigma(i)\sigma(j)}$ for $i, j\in P_{2}$. 

Conversely, it is straightforward to verify that a partition $P_{1}\cup P_{2}$ of $\{1, \cdots, n\}$ and a permutation $\sigma$ of $\{1, \cdots, n\}$ satisfying the conditions can define an algebra automorphism $\varphi$ of $\A$. We will not state the details here. \qed

Next, we solve the isomorphism problem for the family of ``symmetric" multiparameter quantized Weyl algebras. Given two data sets $(n, \overline{q}, \Lambda)$ and $(m, \overline{q}^{\prime}, \Lambda^{\prime})$, we can define two ``symmetric" multiparameter quantized Weyl algebras $\A$ and $\C$. We have the following result.

\begin{thm}
Let $\A$ and $\C$ be two ``symmetric" multiparameter quantized Weyl algebras generated by $x_{1}, y_{1}, \cdots, x_{n}, y_{n}$ and $x_{1}^{\prime}, y_{1}^{\prime}, \cdots, x_{m}^{\prime}, y_{m}^{\prime}$ respectively. Assume that none of the parameters $q_{i}$ and $q_{i}^{\prime}$ is a root of unity, then $\A$ is isomorphic to $\C$ as a $\K-$algebra if and only if $m=n$ and there exist a partition $P_{1}\cup P_{2}$ of $\{1, \cdots, n\}$ and a permutation $\sigma$ of $\{1, 2, \cdots, n\}$ such that
\begin{enumerate}
\item for $i, j\in P_{1}$, we have $ q_{i}=q_{\sigma(i)}^{\prime}, q_{j}=q_{\sigma(j)}^{\prime}, \lambda_{ij}=\lambda_{\sigma(i) \sigma(j)}^{\prime}$. 

\item for $i\in P_{1}$ and $j\in P_{2}$, we have $ q_{i}=q_{\sigma(i)}^{\prime}, q_{j}^{-1}=q_{\sigma(j)}^{\prime},\lambda_{ij}=\lambda_{\sigma(j)\sigma(i)}^{\prime}$.

\item for $i, j\in P_{2}$, we have $ q_{i}^{-1}=q_{\sigma(i)}^{\prime}, q_{j}^{-1}=q_{\sigma(j)}^{\prime},\lambda_{ij}=\lambda_{\sigma(i)\sigma(j)}^{\prime}$.
\end{enumerate}
In particular, the corresponding isomorphism $\varphi \colon \A \longrightarrow \C$ is defined as follows:
\[
\varphi(x_{i})=\alpha_{i} x^{\prime}_{\sigma(i)}, \quad \varphi(y_{i})=\alpha_{i}^{-1} y^{\prime}_{\sigma(i)}
\]
for $i\in P_{1}$ and
\[
\varphi(x_{i})=\alpha_{i} y^{\prime}_{\sigma(i)},\quad \varphi(y_{i})=-q_{i}^{-1} \alpha_{i}^{-1} x^{\prime}_{\sigma(i)}
\]
for $i\in P_{2}$. 
\end{thm}
{\bf Proof:} First of all, the Gelfand-Kirillov dimension of $\A$ is $2n$ and the Gelfand-Kirillov dimension of $\C$ is $2m$. If $\A$ is isomorphic to $\C$ as a $\K-$algebra, then they have the same Gelfand-Kirillov dimension. Thus we have $m=n$. The rest of the proof is to mimic the ones for {\bf Theorem 2.1 and Theorem 2.2} word in word. We will not repeat the details here.
\qed

Let $A$ be any $\K-$algebra. A $\K-$algebra automorphism $h$ of the polynomial extension $A[t]=A\otimes_{\K} \K[t]$ is said to be {\it triangular} if there are a $g\in {\rm Aut}_{\K} (A)$ and $c\in \K^{\ast}$ and $r$ in the center of $A$ such that 
\[
h(w)=ct+r
\]
and
\[
h(a)=g(a)\in A
\]
for any $a\in A$. For more details, we refer the reader to \cite{CPWZ1}.

Denote by $\E=\A\otimes_{\K} \K[t]$ and $\G=\C\otimes_{\K} \K[t]$. Then we have the following result.
\begin{thm}
Let $\varphi\colon \E\longrightarrow \E$ be a $\K-$algebra automorphism of $\E$. Then there exist a partition $P_{1}\cup P_{2}$ of $\{1, \cdots, n\}$ and a permutation $\sigma$ of $\{1, 2, \cdots, n\}$ such that
\begin{enumerate}
\item for $i, j\in P_{1}$, we have $ q_{i}=q_{\sigma(i)}^{\prime}, q_{j}=q_{\sigma(j)}^{\prime}, \lambda_{ij}=\lambda_{\sigma(i) \sigma(j)}^{\prime}$. 

\item for $i\in P_{1}$ and $j\in P_{2}$, we have $ q_{i}=q_{\sigma(i)}^{\prime}, q_{j}^{-1}=q_{\sigma(j)}^{\prime},\lambda_{ij}=\lambda_{\sigma(j)\sigma(i)}^{\prime}$.

\item for $i, j\in P_{2}$, we have $ q_{i}^{-1}=q_{\sigma(i)}^{\prime}, q_{j}^{-1}=q_{\sigma(j)}^{\prime},\lambda_{ij}=\lambda_{\sigma(i)\sigma(j)}^{\prime}$.
\end{enumerate}
In particular, $\varphi$ is defined as follows:
\[
\varphi(x_{i})=\alpha_{i} x_{\sigma(i)}, \quad \varphi(y_{i})=\alpha_{i}^{-1} y_{\sigma(i)}
\]
for $i\in P_{1}$ and
\[
\varphi(x_{i})=\alpha_{i} y_{\sigma(i)},\quad \varphi(y_{i})=-q_{i}^{-1} \alpha_{i}^{-1} x_{\sigma(i)}
\]
for $i\in P_{2}$ and
\[
\varphi(t)=at+b
\]
for some $\alpha_{i}, a\in \K^{\ast}$ and $b\in \K$.
\end{thm}
{\bf Proof:} Note that the localization of $\E$ with respect to the Ore set $\K[t]-\{0\}$ is isomorphic to the algebra $\mathcal{A}_{n}^{\overline{q}, \Lambda}(\K(t))$, which is a ``symmetric" multiparameter quantized Weyl algebra over the base field $\K(t)$. Since none of $q_{i}$ is a root of unity, $\mathcal{A}_{n}^{\overline{q}, \Lambda}(\K(t))$ is a simple $\K(t)-$algebra. Thus the center of $\E$ is $\K[t]$. Since the center of $\E$ is preserved by $\varphi$, we can restrict $\varphi$ to a $\K-$algebra automorphism of $\K[t]$. Thus we have $\varphi(t)=at+b$ for some $a\in K^{\ast}$ and $b\in K$.

Let us denote the restriction of $\varphi$ to $\K[t]$ by $\varphi_{0}$. Note that the inverse of $\varphi_{0}$ can be extended to a $\K-$algebra $\psi$ of $\E$ as follows:
\[
\psi(t)=\varphi_{0}^{-1}(t),\quad \psi(w)=w
\]
for any $w\in \A$. Now the composition $\psi\circ \varphi$ is actually a $\K[t]-$algebra automorphism of $\E$. Thus we can extend $\psi\circ \varphi$ to a $\K(t)-$algebra automorphism of $\mathcal{A}_{n}^{\overline{q},\Lambda}(\K(t))$ and we still denote the extension by $\psi\circ \varphi$. By {\bf Theorem 2.2}, there exist a partition $P_{1}\cup P_{2}$ of $\{1, \cdots, n\}$ and a permutation $\sigma$ of $\{1, 2, \cdots, n\}$ such that
\[
\psi\circ \varphi(x_{i})=\alpha_{i} x_{\sigma(i)}, \quad \psi \circ \varphi(y_{i})=\alpha_{i}^{-1} y_{\sigma(i)}
\]
with $\alpha_{i}\in \K(t)^{\ast}$ for $i\in P_{1}$ and
\[
\psi\circ \varphi(x_{i})=\alpha_{i} y_{\sigma(i)},\quad \psi\circ \varphi(y_{i})=-q_{i}^{-1} \alpha_{i}^{-1} x_{\sigma(i)}
\]
with $\alpha_{i} \in \K(t)^{\ast}$ for $i\in P_{2}$. Since $\psi\circ \varphi (x_{i})\in \A$, we have that $\alpha_{i}\in \K^{\ast}$. Thus, we are done with the proof.
\qed

\begin{thm}
Any $\K-$algebra automorphism of $\E$ is triangular. In particular, we have the following:
\[
{\rm Aut}_{\K}(\A[t])=\left(\begin{array}{lr} 
{\rm Aut}_{\K}(\A) & \K\\
0 & \K^{\ast}
\end{array}
\right).
\]
As a result, if $\Z \subset \K$, then we have ${\rm LNDer} (\A)=\{0\}$, where ${\rm LNDer} (\A)$ is the set of all locally nilpotent derivations of $\A$. 
\end{thm}
{\bf Proof:} By {\bf Theorem 2.5}, every $\K-$algebra automorphism $\varphi$ of $\E$ can be restricted to $\A$ and $\varphi(t)=at+b$ for some $a\in \K^{\ast}$ and $b\in K$. Thus $\varphi$ is triangular. As a result, we have the following:
\[
{\rm Aut}_{\K}(\E)=\left(\begin{array}{lr} 
{\rm Aut}_{\K}(\A) & \K\\
0 & \K^{\ast}
\end{array}
\right).
\]

Let $\partial$ be a locally nilpotent derivation of $\A$. Since $\Z\subset \K$, we can define a $\K(t)-$algebra automorphism $\varphi$ of $\mathcal{A}_{n}^{\overline{q}, \Lambda}(\K(t))$ as follows:
\[
\varphi(t)=t, \quad\varphi(x)=\sum_{i=0}^{\infty} \frac{t^{i}}{i!} \partial^{i}(x)
\]
for any $x\in \A$. If $\partial\neq 0$, then $\partial(x_{i})\neq 0$ for some $i$ or $\partial (y_{j})\neq 0$ for some $j$. Then we have
\[
\varphi(x_{i})=x_{i}+t\partial(x_{i})+\text{other terms\, involving higher powers of }\, t
\]
or 
\[
\varphi(y_{j})=y_{j}+t\partial(y_{j})+\text{other terms involving higher powers of }\, t.
\]
This is a contradiction to the description of the $\K(t)-$automorphisms of $\mathcal{A}_{n}^{\overline{q}, \Lambda}(\K(t))$. Therefore, we have $\partial=0$. So we have completed the proof.
\qed

\begin{thm}
The algebra $\E$ is isomorphic to $\G$ if and only if $m=n$, and there exist a partition $P_{1}\cup P_{2}$ of $\{1, \cdots, n\}$ and a permutation $\sigma$ of $\{1, 2, \cdots, n\}$ such that
\begin{enumerate}
\item for $i, j\in P_{1}$, we have $ q_{i}=q_{\sigma(i)}^{\prime}, q_{j}=q_{\sigma(j)}^{\prime}, \lambda_{ij}=\lambda_{\sigma(i) \sigma(j)}^{\prime}$. 

\item for $i\in P_{1}$ and $j\in P_{2}$, we have $ q_{i}=q_{\sigma(i)}^{\prime}, q_{j}^{-1}=q_{\sigma(j)}^{\prime},\lambda_{ij}=\lambda_{\sigma(j)\sigma(i)}^{\prime}$.

\item for $i, j\in P_{2}$, we have $ q_{i}^{-1}=q_{\sigma(i)}^{\prime}, q_{j}^{-1}=q_{\sigma(j)}^{\prime},\lambda_{ij}=\lambda_{\sigma(i)\sigma(j)}^{\prime}$.
\end{enumerate}
In particular, the corresponding isomorphism $\varphi \colon \E \longrightarrow \G$ is defined as follows:
\[
\varphi(x_{i})=\alpha_{i} x^{\prime}_{\sigma(i)}, \quad \varphi(y_{i})=\alpha_{i}^{-1} y^{\prime}_{\sigma(i)}
\]
for $i\in P_{1}$ and
\[
\varphi(x_{i})=\alpha_{i} y^{\prime}_{\sigma(i)},\quad \varphi(y_{i})=-q_{i}^{-1} \alpha_{i}^{-1} x^{\prime}_{\sigma(i)}
\]
for $i\in P_{2}$ and
\[
\varphi(t)=at+b
\]
for some $a\in \K^{\ast}$ and $b\in \K$.
\end{thm}

{\bf Proof:} Note that $\A$ is cancellative. Then $\E$ is isomorphic to $\G$ if and only if $\A$ is isomorphic to $\C$. Now the result follows from {\bf Theorem 2.3}.

\qed

Let $\B$ and $\D$ be two Maltsiniotis multiparameter quantized Weyl algebras generated by $x_{1}, y_{1}, \cdots, x_{n}, y_{n}$ and $x_{1}^{\prime}, y_{1}^{\prime}, \cdots, x_{m}^{\prime}, y_{m}^{\prime}$ respectively. Let us set $\F=\B\otimes_{\K} \K[t]$ and $\H=\D\otimes_{\K} \K[t]$. We have the following result for the automorphism group of the polynomial extension of the Maltsiniotis multiparameter quantized Weyl algebra.

\begin{thm}
Assume that none of $q_{i}, i=1, 2, \cdots, n$ is a root of unity. Then we have the following results.
\begin{itemize}
\item If $\varphi $ is a $\K-$algebra automorphism of $\F$, then there exists a tuple $(\mu_{1}, \cdots, \mu_{n}, a, b)\in (\K^{\ast})^{n+1}\times \K$ such that $ \varphi(t)=at+b$ and $\varphi(x_{i})=\mu_{i}x_{i}$, $\varphi(y_{i})=\mu_{i}^{-1}y_{i}$ for each $1\leq i\leq n$.\\

\item For any $a \in \K^{\ast}$ and $b\in \K$, one can define a $\K-$algebra automorphism $\varphi_{a,b}$ of $\F$ as follows:
\[
\varphi_{a,b}(t)=at+b,\quad \varphi_{a,b}(x_{i})=x_{i},\quad \varphi_{a,b}(y_{i})=y_{i}
\]
for each $1\leq i\leq n$.

\item We have that ${\rm Aut}_{\K}(\F)\cong (\K^{\ast})^{n}\rtimes G$ where $G$ is the subgroup of ${\rm Aut}_{\K}(\F)$ consisting of these $\varphi_{a,b}$.
\end{itemize}
\end{thm}

{\bf Proof:} The proof is similar to the one for {\bf Theorem 2.4} and we will not repeat the details here.\qed

Moreover, we have the following theorem.
\begin{thm}
Any $\K-$algebra automorphism of $\F$ is triangular. In particular, we have the following:
\[
{\rm Aut}_{\K}(\F)=\left(\begin{array}{lr} 
{\rm Aut}_{\K}(\B) & \K\\
0 & \K^{\ast}
\end{array}
\right).
\]
As a result, if $\Z \subset \K$, then we have ${\rm LNDer} (\B)=\{0\}$, where ${\rm LNDer} (\B)$ is the set of all locally nilpotent derivations of $\B$. 
\end{thm}
{\bf Proof:} The proof is the same as the one used for {\bf Theorem 2.5}. We skip the details.\qed

Next, we solve the isomorphism classification problem for the family of polynomial extensions $\{\E\}$ based on a result due to Goodearl and Hartwig \cite{GH}. For convenience, we will follow the notation used in \cite{GH}.
\begin{thm}
Assume none of $q_{1}, \cdots, q_{n}$ and $q_{1}^{\prime}, \cdots, q_{m}^{\prime}$ is a root of unity. Then $\F$ is isomorphic to $\H$ as a $\K-$algebra if and only if 
\begin{itemize}
\item $m=n$;
\item There exists a sign vector $\varepsilon \in \{-1,1\}^{n}$ such that we have $q_{i}^{\prime}=q_{i}^{\varepsilon_{i}}, \forall 1\leq i \leq n$ and
\begin{displaymath}
\gamma_{ij}^{\prime}=
\left \{ 
\begin{array}{ll}
\gamma_{ij} & \, \text{if}\, (\varepsilon_{i}, \varepsilon_{j})=(1,1)\\
\gamma_{ji} & \, \text{if}\, (\varepsilon_{i}, \varepsilon_{j})=(-1,1) \\
q_{i}^{-1}\gamma_{ji} & \, \text{if}\, (\varepsilon_{i}, \varepsilon_{j})=(1,-1)\\
q_{i}\gamma_{ij}& \, \text{if}\, (\varepsilon_{i}, \varepsilon_{j})=(-1, -1)
\end{array}
\right.\,\forall 1\leq i\leq j\leq n.
\end{displaymath}
\end{itemize}
If the above conditions are satisfied, then for any $\mu \in (\K^{\ast})^{n}$ and $\varepsilon \in \{\pm 1\}^{n}$ and $a\in \K^{\ast}$ and $b\in \K$, one can define a unique $\K-$algebra isomorphism $\varphi_{\mu, \varepsilon, a, b}\colon \E\longrightarrow \H$
by
\[
\varphi_{\mu, \varepsilon, a, b}(t)=at+b
\]
and
\[
\varphi_{\mu, \varepsilon, a, b}(x_{i})=\left \{
\begin{array}{ll} \mu_{i} x_{i}^{\prime},\quad \varepsilon_{i}=1\\
\mu_{i}y_{i}^{\prime},\quad \varepsilon_{i}=-1
\end{array}\right.
\]
and
\[
\varphi_{\mu, \varepsilon, a, b}(y_{i})=\left \{
\begin{array}{ll} \lambda_{i}\mu_{i}^{-1} y_{i}^{\prime},\quad\quad \varepsilon_{i}=1\\
-\lambda_{i}\mu_{i}^{-1} x_{i}^{\prime},\quad \varepsilon_{i}=-1
\end{array} \right.
\]
for each $i\in \{1, \cdots, n\}$ with $\lambda=\lambda(\varepsilon)\in (\K^{\ast})^{n}$ recursively defined as follows:
\[
\lambda_{0}=1,\quad \lambda_{i} =q^{(\varepsilon_{i}-1)/2}\lambda_{i-1}.
\]

\end{thm}
{\bf Proof:} Note that the algebra $\B$ is cancellative. Thus we have that $\F\cong \H$ if and only if $\B\cong \D$. Now the result follows from {\bf Theorem 5.1} in \cite{GH}.\qed

\section{The Quantum Dixmier Conjecture for $(\A)_{\mathcal{Z}}$}
In this section, we prove a quantum analogue of the Dixmier conjecture for the simple localization $(\A)_{\mathcal{Z}}$ of $\A$ under the condition that $q_{1}^{i_{1}}q_{2}^{i_{2}}\cdots q_{n}^{i_{n}}=1$ implies $i_{1}=i_{2}=\cdots =i_{n}=0$. In particular, we will show that each $\K-$algebra endomorphism of $(\A)_{\mathcal{Z}}$ is an algebra automorphism. Furthermore, we will describe the automorphism group for $(\A)_{\mathcal{Z}}$. 

\begin{thm} Suppose that $q_{1}^{i_{1}}q_{2}^{i_{2}}\cdots q_{n}^{i_{n}}=1$ implies $i_{1}=i_{2}=\cdots =i_{n}=0$. Then every $\K-$algebra endomorphism $\varphi$ of $(\A)_{\mathcal{Z}}$ is a $\K-$algebra automorphism.\\

\end{thm}

{\bf Proof:} Note that any invertible element of $(\A)_{\mathcal{Z}}$ is of the form $\alpha z_{1}^{a_{1}}\cdots z_{n}^{a_{n}}$, where $\alpha \in \K^{\ast}$ and $a_{1}, \cdots, a_{n}\in \Z$. Let $\varphi$ be a $\K-$algebra endomorphism of $(\A)_{\mathcal{Z}}$, then $\varphi$ sends invertible elements of $(\A)_{\mathcal{Z}}$ to its invertible elements. Therefore, for $i=1, \cdots, n$, we have the following
\[
\sigma(z_{i})=\lambda_{i} z_{1}^{a_{i1}}\cdots z_{n}^{a_{in}}
\]
where $\lambda_{i}\in \K^{\ast}$ and $a_{i1}, \cdots, a_{in}\in \Z$. Since $y_{i}x_{i}=\frac{z_{i}-1}{q_{i}-1}, x_{i}y_{i}=\frac{q_{i}z_{i}-1}{q_{i}-1}$ and $z_{i}x_{i}=q_{i}^{-1}x_{i}z_{i}, z_{i}y_{i}=q_{i}y_{i}z_{i}$ for $i=1, \cdots, n$, we have that $\sigma(z_{i})\notin \K^{\ast}$ for $i=1, \cdots, n$.

Note that we can realize $(\A)_{\mathcal{Z}}$ as a generalized Weyl algebra over the base algebra $\K[z_{1}^{\pm 1}, \cdots, z_{n}^{\pm 1}]$ with the automorphisms $\sigma_{i}$ of $\K[z_{1}^{\pm 1}, \cdots, z_{n}^{\pm 1}]$ defined by $\sigma_{i}(z_{j})=q_{i}^{\delta_{ij}} z_{j}$. Thus for $i=1, \cdots, n$, we can write the image of $x_{i}$ under $\varphi$ as follows:
\[
\sigma(x_{i})=\sum_{\overline{k}^{i}} f_{\overline{k}^{i}} x_{1}^{k_{1}^{i}}\cdots x_{n}^{k_{n}^{i}}+ \sum_{\overline{l}^{i}}g_{\overline{l}^{i}} y_{1}^{l_{1}^{i}}\cdots y_{n}^{l_{n}^{i}}
\]
where $\overline{k}^{i}=(k^{i}_{1}, \cdots, k^{i}_{n})\in (\Z_{\geq 0})^{n}$, and $\overline{l}^{i}=(l^{i}_{1}, \cdots, l^{i}_{n}) \in (\Z_{\geq 0})^{n}$, and $f_{\overline{k}^{i}}\in \K[z_{1}^{\pm 1}, \cdots, z_{n}^{\pm 1}]$, and $g_{\overline{l}^{i}}\in \K[z_{1}^{\pm 1}, \cdots, z_{n}^{\pm 1}]$. 

Since $z_{i}x_{i}=q_{i}^{-1}x_{i}z_{i}$, we have that
\[
q_{1}^{-a_{i1} k^{i}_{1}}q_{2}^{-a_{i2}k^{i}_{2}}\cdots q_{n}^{-a_{in}k^{i}_{n}}=q_{i}^{-1}
\]
for $f_{\overline{k}^{i}}\neq 0$ and
\[
q_{1}^{a_{i1}l^{i}_{1}}\cdots q_{n}^{a_{in}l^{i}_{n}}=q_{i}^{-1}
\]
for $g_{\overline{l}^{i}}\neq 0$. Since $q_{1}^{i_{1}}q_{2}^{i_{2}}\cdots q_{n}^{i_{n}}=1$ implies $i_{1}=\cdots =i_{n}=0$, we have that $a_{ii}k^{i}_{i}=1$ for $f_{\overline{k}^{i}}\neq 0$ and $a_{ii}l^{i}_{i}=-1$ for $g_{\overline{l}}\neq 0$. Since $\overline{k}^{i}\in (\Z_{\geq 0})^{n}$ and $\overline{l}^{i}\in (\Z_{\geq 0})^{n}$, we cannot have both $f_{\overline{k}^{i}}\neq 0$ for some $\overline{k}^{i}$ and $g_{\overline{l}^{i}}\neq 0$ for some $\overline{l}^{i}$ spontaneously. Thus for $i=1, \cdots, n$, we have either
\[
\sigma(x_{i})=\sum_{\overline{k}^{i}} f_{\overline{k^{i}}} x_{1}^{k^{i}_{1}}\cdots x_{n}^{k^{i}_{n}}
\]
or
\[
\sigma(x_{i})=\sum_{\overline{l}^{i}}g_{\overline{l^{i}}} y_{1}^{l^{i}_{1}}\cdots y_{n}^{l^{i}_{n}}.
\]
Since $z_{j}x_{i}=x_{i}z_{j}$ for $i\neq j$, we have $a_{ji}k^{i}_{i}=0$ or $a_{ji}l^{i}_{i}=0$ for $j\neq i$. Thus we have that $a_{ji}=0$ for $j\neq i$ and $a_{ii}=\pm 1$. As a result, we have 
\[
\varphi(z_{i})=\lambda_{i} z_{i},\quad \text{or}\quad \varphi(z_{i})=\lambda_{i} z_{i}^{-1}
\]
for $i=1, \cdots, n$. 

Suppose that we have 
\[
\sigma(x_{i})=\sum_{\overline{k}^{i}} f_{\overline{k}^{i}} x_{1}^{k^{i}_{1}}\cdots x_{n}^{k^{i}_{n}}.
\]

Since $k^{i}_{i}\geq 0$, we have that $a_{ii}=1$. So we can conclude that $\varphi(z_{i})=\lambda_{i} z_{i}$. Since $z_{i}x_{j}=x_{j}z_{i}$ and $z_{i}y_{j}=y_{j}z_{j}$ for $i\neq j$, we have that $\varphi(x_{i})=f_{i} x_{i}$, and $\varphi(y_{i})=g_{i} y_{i}$ for some $f_{i}, g_{i} \in \K[z_{1}^{\pm 1}, \cdots, z_{n}^{\pm 1}]$. Since we have $x_{i}y_{i}=\frac{q_{i}z_{i}-1}{q_{i}-1}$, we have the following:
\[
f_{i}x_{i}g_{i} y_{i}=\frac{\lambda_{i}q_{i}z_{i}-1}{q_{i}-1}.
\]
Thus we have 
\[
f_{i}x_{i}g_{i} y_{i}=f_{i}g_{i}^{\prime} x_{i}y_{i}=f_{i}g_{i}^{\prime} \frac{q_{i}z_{i}-1}{q_{i}-1}=\frac{\lambda_{i} q_{i}z_{i}-1}{q_{i}-1}
\]
for some $g_{i}^{\prime}\in \K[z_{1}^{\pm 1}, \cdots, z_{n}^{\pm 1}]$. As a result, we have $f_{i}g_{i}^{\prime}=1$ and $\lambda_{i}=1$. Note that $g_{i}^{\prime}=\sigma_{i}(g_{i})$. Therefore, we have that $f_{i}=\mu_{i}z_{1}^{b_{i1}}\cdots z_{n}^{b_{in}}$ and $g_{i}=\mu_{i}^{-1}q_{i}^{b_{ii}} z_{1}^{-b_{i1}}\cdots z_{n}^{-b_{in}}$ for some $\mu_{i} \in \K^{\ast}$ and $b_{i1}, \cdots, b_{in}\in \Z$.

Suppose that we have 
\[
\sigma(x_{i})=\sum_{\overline{l}^{i}}g_{\overline{l}^{i}} y_{1}^{l^{i}_{1}}\cdots y_{n}^{l^{i}_{n}}.
\]
Then we have $a_{ji}l^{i}_{i}=0$ for $j\neq i$ and $a_{ii}l^{i}_{i}=-1$. Thus we have $a_{ji}=0$ for $j\neq i$. Since $l^{i}_{i}\in \Z_{\geq 0}$, we have $a_{ii}=-1$. As a result, we have $\varphi(z_{i})=\lambda_{i} z_{i}^{-1}$ for some $\lambda_{i}\in \K^{\ast}$. Since $z_{j}x_{i}=x_{i}z_{j}$ for $i\neq j$, we have that $\varphi(x_{i})=g_{i} y_{i}$ and $\varphi(y_{i})=f_{i} x_{i}$ for some $f_{i}, g_{i}\in \K[z_{1}^{\pm 1}, \cdots, z_{n}^{\pm 1}]$. Once again, since $x_{i}y_{i}=\frac{q_{i}z_{i}-1}{q_{i}-1}$, we have following:
\[
g_{i}y_{i}f_{i}x_{i}=\frac{\lambda_{i}q_{i}z_{i}^{-1}-1}{q_{i}-1}.
\]
Since $y_{i}x_{i}=\frac{z_{i}-1}{q_{i}-1}$, we have that
\[
g_{i}y_{i}f_{i}x_{i}=g_{i}f_{i}^{\prime} y_{i}x_{i}=g_{i}f_{i}^{\prime}\frac{z_{i}-1}{q_{i}-1}=\frac{\lambda_{i}q_{i}z_{i}^{-1}-1}{q_{i}-1}
\]
for some $f_{i}^{\prime}\in \K[z_{1}^{\pm 1}, \cdots, z_{n}^{\pm 1}]$. As a result, we have $g_{i}f_{i}^{\prime}z=-1$ and $\lambda_{i}=q_{i}^{-1}$. Since $f_{i}^{\prime}=\sigma_{i}^{-1}(f_{i})$, we have that $f_{i}=\mu_{i}z_{1}^{b_{i1}}\cdots z_{i}^{b_{ii}}\cdots z_{n}^{b_{in}}$ and $h_{i}=-\mu_{i}^{-1}q_{i}^{-b_{ii}-1} z_{1}^{-b_{i1}}\cdots z_{i}^{-b_{ii}-1}\cdots z_{n}^{-b_{in}}$ for some $\mu_{i}\in \K^{\ast}$ and $b_{i1}, \cdots, b_{in}\in \Z$.

Let $\varphi$ be any $\K-$algebra endomorphism of $(\A)_{\mathcal{Z}}$, then for each $i=1, \cdots, n$, we have either $\varphi(z_{i})=\lambda_{i} z_{i}$ or $\varphi(z_{i})=\lambda_{i}z_{i}^{-1}$ for some $\lambda_{i}\in \K^{\ast}$. We will say that a $\K-$algebra endomorphism $\varphi$ of $(\A)_{\mathcal{Z}}$ is of positive type if 
\[
\varphi(z_{i})=\lambda_{i}z_{i}
\]
for $i=1, \cdots, n$. Even though $\varphi$ may not be of positive type, $\varphi\circ \varphi$ is a $\K-$algebra endomorphism of $(\A)_{\mathcal{Z}}$ of positive type. Next, we will show that each $\K-$algebra endomorphism $\varphi$ of $(\A)_{\mathcal{Z}}$ of positive type is an automorphism, which implies that any $\K-$algebra endomorphism $\varphi$ of $(\A)_{\mathcal{Z}}$ is an automorphism.

Let $\varphi$ be a $\K-$algebra endomorphism of $(\A)_{\mathcal{Z}}$ of positive type, then for each $i=1, \cdots, n$, we have the following:
\[
\varphi(x_{i})=\mu_{i}z_{1}^{b_{i1}}\cdots z_{n}^{b_{in}}x_{i},\quad \varphi(y_{i})=\mu_{i}^{-1}q_{i}^{b_{ii}} z_{1}^{-b_{i1}}\cdots z_{n}^{-b_{in}}y_{i}, \quad \varphi(z_{i})=z_{i}
\]
for some $\mu_{i} \in \K^{\ast}$ and $b_{i1}, \cdots, b_{in}\in \Z$.

Note that $q_{1}^{i_{1}}\cdots q_{n}^{i_{n}}=1$ implies that $i_{1}=\cdots =i_{n}=0$ and 
\[
x_{i}x_{j}=\lambda_{ij} x_{j}x_{i},\quad x_{i}y_{j}=\lambda_{ji} y_{j}x_{i}, \quad y_{i}y_{j}=\lambda_{ij} y_{j}y_{i}
\]
for $i\neq j$. After applying $\varphi$ to both sides of the above commuting relations, we can conclude that 
\[
\varphi(x_{i})=\mu_{i}z_{i}^{b_{ii}} x_{i}, \quad \varphi(y_{i})=\mu_{i}^{-1}q_{i}^{b_{ii}} z_{i}^{-b_{ii}} y_{i},\quad \varphi(z_{i})=z_{i}
\]
for $i=1, \cdots, n$. It is obvious that $\varphi$ has an inverse. Thus $\varphi$ is a $\K-$algebra automorphism of $(\A)_{\mathcal{Z}}$. So we have proved that every $\K-$algebra endomorphism of $(\A)_{\mathcal{Z}}$ is an automorphism.
\qed

\begin{thm} Maintain the condition that $q_{1}^{i_{1}}\cdots q_{n}^{i_{n}}=1$ implies $i_{1}=\cdots =i_{n}=0$ and let $\varphi$ be a $\K-$algebra automorphism of $(\A)_{\mathcal{Z}}$. Then there exists a partition $P_{1}\cup P_{2}$ of the set $\{1, \cdots, n\}$ with $P_{1}=\{l_{1}, \dots, l_{k}\}$ and $P_{2}=\{l_{k+1}, \cdots, k_{n}\}$ such that 
\begin{enumerate}
\item For any $l_{i}\in P_{1}$ and $l_{j}\in P_{2}$ , we have that
\[
\varphi(x_{l_{i}})=\alpha_{l_{i}} z_{l_{k+1}}^{a_{i(k+1)}} \cdots z_{l_{n}}^{a_{in}} z_{l_{i}}^{c_{i}}x_{l_{i}}
\]
and
\[
\varphi(y_{l_{i}})=\beta_{l_{i}} z_{l_{k+1}}^{-a_{i(k+1)}} \cdots z_{l_{n}}^{-a_{in}}z_{l_{i}}^{-c_{i}}y_{l_{i}}
\]
and
\[
\varphi(x_{l_{j}})=\alpha_{l_{j}} z_{l_{1}}^{b_{j 1}}\cdots z_{l_{k}}^{b_{jk}} z_{l_{j}}^{c_{j}}y_{l_{j}}
\]
and
\[
\varphi(y_{j})=\beta_{l_{j}} z_{l_{1}}^{-b_{j1}} \cdots z_{l_{k}}^{-b_{jk}} z_{l_{j}}^{-c_{j}-1}x_{l_{j}}
\]
for some $a_{i(k+1)}, \cdots, a_{in}, b_{j1}, \cdots, b_{jk}\in \Z$, and $c_{i}, c_{j} \in \mathbb{Z}$, and $\alpha_{l_{i}}, \beta_{l_{i}}, \alpha_{l_{j}}, \beta_{l_{j}}\in \mathbb{K}^{\ast}$ such that 
\[
\alpha_{l_{i}}\beta_{l_{i}}=q_{l_{i}}^{c_{i}},\quad
\alpha_{l_{j}}\beta_{l_{j}}=-q_{l_{j}}^{-c_{j}-1},\quad
q_{l_{i}}^{b_{ji}} q_{l_{j}}^{a_{ij}}=\lambda_{l_{i}l_{j}}^{2}.
\]

\item For any given partition of $P_{1}\cup P_{2}$ of the set $\{1, \cdots, n\}$ satisfying the condition, let $\varphi$ be a $\mathbb{K}-$linear mapping defined as above. Then $\varphi$ can be extended to a $\K-$algebra automorphism of $(\A)_{\mathcal{Z}}$.
\end{enumerate}
\end{thm}
{\bf Proof:} Let $\varphi$ be any $\K-$algebra automorphism of $(\A)_{\mathcal{Z}}$. Then for $i=1, \cdots, n$, we have either

\noindent
{\bf Case 1:} $\varphi(x_{i})=\mu_{i}z_{1}^{b_{i1}}\cdots z_{n}^{b_{in}}x_{i}$, and $\varphi(y_{i})=\mu_{i}^{-1}q_{i}^{b_{ii}} z_{1}^{-b_{i1}}\cdots z_{n}^{-b_{in}}y_{i}$, and $\varphi(z_{i})=z_{i}$;
or 

\noindent
{\bf Case 2:} $\varphi(z_{i})=q_{i}^{-1}z_{i}^{-1}$, and $\varphi(x_{i})=\mu_{i}z_{1}^{b_{i1}}\cdots z_{i}^{b_{ii}}\cdots z_{n}^{b_{in}}y_{i}$, 
and $\varphi(y_{i})=-\mu_{i}^{-1}q_{i}^{-b_{ii}-1} z_{1}^{-b_{i1}}\cdots z_{i}^{-b_{ii}-1}\cdots z_{n}^{-b_{in}}x_{i}$.

We can set $P_{1}$ to be the set consisting of all $i\in \{1, \cdots, n\}$ such that the first case is true and $P_{2}$ to be the set consisting of all $i\in \{1, \cdots, n\}$ such that the second case is true. It is obvious that $P_{1}$ and $P_{2}$ are disjoint and $P_{1}\cup P_{2}$ is a partition of $\{1, \cdots, n\}$. The rest is verified by applying the automorphism $\varphi$ to both sides of the following relations: 
\[
x_{i}x_{j}=\lambda_{ij} x_{j} x_{i}, \quad x_{i}y_{j}=\lambda_{ji} y_{j}x_{i},\quad y_{i}y_{j}=\lambda_{ij} y_{j}y_{i}
\]
for any $i\neq j$. 

Conversely, it is straightforward to verify that the mapping $\varphi$ defined by any given partition $P_{1}\cup P_{2}$ satisfying the condition can be extended to a $\K-$algebra automorphism of $(\A)_{\mathcal{Z}}$. 
\qed

As an immediate consequence of {\bf Theorem 3.2}, we have the following corollary on the automorphism group of $(\A)_{\mathcal{Z}}$ in some special cases.

\begin{cor} Let us further assume that $q_{i}^{k_{i}}q_{j}^{k_{j}}=\lambda_{ij}^{2}$ implies $k_{i}=k_{j}=0$ for any $i\neq j$, then we have the following result.
\begin{enumerate}
\item Let $\varphi$ be a $\K-$algebra automorphism $\varphi$ of $(\A)_{\mathcal{Z}}$. For each $i=1, \cdots, n$, we have {\bf either}
\[
\varphi(x_{i})=\alpha_{i} z_{i}^{c_{i}} x_{i},\quad\quad \varphi(y_{i})=\beta_{i} z_{i}^{-c_{i}}y_{i}
\]
for some $c_{i}\in \mathbb{Z}$ and $\alpha_{i}, \beta_{i}\in \mathbb{K}^{\ast}$ such that $\alpha_{i}\beta_{i}=q_{i}^{c_{i}}$, {\bf or}
\[
\varphi(x_{i})=\alpha_{i} z_{i}^{c_{i}} y_{i}, \quad\quad \varphi(y_{i})=\beta_{i} z_{i}^{-c_{i}-1}x_{i}
\]
for some $c_{i}\in \mathbb{Z}$ and $\alpha_{i}, \beta_{i}\in \mathbb{K}^{\ast}$ such that $\alpha_{i}\beta_{i}=-q_{i}^{-c_{i}-1}$.\\

\item Let $\varphi$ be a $\mathbb{K}-$linear mapping defined as above. Then $\varphi$ can be extended to a $\K-$algebra automorphism of $(\A)_{\mathcal{Z}}$.
\end{enumerate}
\end{cor}
\qed

\noindent
{\bf Acknowledgements:} The author would like to thank James Zhang for sharing an earlier version of the preprint \cite{BZ}. Part of the results were announced during the special session on ``New Developments in Non-commutative Algebra", the AMS Spring Western Sectional Meeting at University of Nevada, Las Vegas, Las Vegas, NV, April 18-19, 2015. The author would like to thank the organizers for the hospitality.


\begin{thebibliography}{99999999}
\frenchspacing

\bibitem{AJ} M. Akhavizadegan and D. A. Jordan, Prime ideals of quantized Weyl algebras, {\it Glasgow Math. J.}, {\bf 38} (1996), no. 3, 283--297.

\bibitem{AD} J. Alev and F. Dumas, Rigidit\'{e} des plongements des quotients primitifs minimaux de $U_{q} (\mathfrak{s}\mathfrak{l} ( 2 ))$ dans l’alg\`{e}bre quantique de Weyl-Hayashi, {\it Nagoya Math. J.}, {\bf 143} (1996),119--146.

\bibitem{Backelin} E. Backelin, Endomorphisms of quantized Weyl algebras, {\it Lett. Math. Phys.}, {\bf 97} (2011), no.
3, 317--338.

\bibitem{Bavula} V. V. Bavula, Generalized Weyl algebras and their representations, {\it St. Petersburg
Math. J.}, {\bf 4} (1993), no. 1, 71--92.

\bibitem{BJ} V. V. Bavula and D. A. Jordan, Isomorphism problems and groups of automorphisms for
generalized Weyl algebras, {\it Trans. Amer. Math. Soc.}, {\bf 353} (2001), no. 2, 769--794.


\bibitem{BZ} J. P. Bell and J. J. Zhang, Zariski cancellation problem for noncommutative algebras, preprint, arXiv:1601.04625.


\bibitem{BK} A. Belov-Kanel and M. Kontsevich, The Jacobian conjecture is stably equivalent to the Dixmier conjecture, {\it Mosc. Math. J.}, {\bf 7} (2007), no. 2, 209--218.


\bibitem{CPWZ1} S. Ceken, J. H. Palmieri, Y. H. Wang and J. J. Zhang, The discriminant controls automorphism groups of noncommutative algebras, {\it Adv. Math.}, {\bf 269} (2015), 551--584.

\bibitem{CPWZ2} S. Ceken, J. H. Palmieri, Y. H. Wang and J. J. Zhang, The discriminant criterion and automorphism groups of quantized algebras, {\it Adv. Math.} {\bf 286} (2016), 754--801

\bibitem{CKZ} K. Chan, A. Young and J. Zhang, Discriminant formulas and applications, preprint, arXiv:1503.06327.


\bibitem{Dixmier} J. Dixmier, Sur les alg\'{e}bres de Weyl, {\it Bull. Soc. Math. France}, {\bf 96} (1968), 209--242.


\bibitem{FKK} H. Fujita, E. Kirkman and J. Kuzmanovich, Global and Krull dimensions of quantum Weyl algebras, {\it J. Algebra}, {\bf 216} (1999), 405--416.

\bibitem{GKM} J. Gaddis, E. Kirkman and W. F. Moore, On the discriminant of twisted tensor products, preprint, arXiv:1606.03105.

\bibitem{GZ} A. Giaquinto and J. J. Zhang, Quantum Weyl algebras, {\it J. Algebra}, {\bf 176} (1995), 861--881.

\bibitem{Ginzburg} V. Ginzburg, Calabi-Yau algebras, preprint, arXiv:0612139v3.


\bibitem{GK} J. Gomez-Torrecillas and L. EL Kaotit, The group of automorphisms of the coordinate ring of quantum symplectic space, {\it Beitr. Algebra Geom.}, {\bf 43} (2002), 597--601.

\bibitem{Goodearl} K. R. Goodearl, Prime ideals in skew polynomial rings and quantized Weyl algebras, {\it J. Algebra} {\bf 150} (1992), no. 2, 324--377. 

\bibitem{GH} K. R. Goodearl and J. T. Hartwig, The isomorphism problem for multiparameter quantized Weyl algebras, {\it S\~{a}o Paulo J. Math. Sci.}, {\bf 9} (2015), 53--61.


\bibitem{GY} K. R. Goodearl and M. T. Yakimov, Unipotent and Nakayama automorphisms of quantum nilpotent algebras, preprint, arXiv:1311.0278.


\bibitem{Jordan} D. A. Jordan, A simple localization of the quantized Weyl algebra, {\it J. Algebra}, {\bf 174} (1995), no. 1, 267--281.

\bibitem{Keller} O. H. Keller, Ganze Cremona-Transformationen, {\it Monatsh. Math. Phys.}, {\bf 47} (1939), no. 1,
299--306.

\bibitem{KL} A. P. Kitchin and S. Launois, Endomorphisms of quantum generalized Weyl algebras, {\it Lett. Math. Phys.}, {\bf 104} (2014), 837--848.

\bibitem{KL1}A. P. Kitchin and S. Launois, On the automorphisms of quantum Weyl algebras, preprint, arXiv:1511.01775.

\bibitem{LY} J. Levitt and M. Yakimov, Quantized Weyl algebras at roots of unity, preprint, arXiv:1606.02121.

\bibitem{LWW} L.Y. Liu, S. Q. Wang and Q. S. Wu, Twisted Calabi-Yau property of Ore extensions, {\it J. Noncommut. Geom.}, {\bf 8} (2014), no. 2, 587--609.

\bibitem{Maltsiniotis} G. Maltsiniotis, Groupes quantiques et structures différentielles, {\it C. R. Acad. Sci. Paris S\'{e}r. I Math.}, {\bf 311} (1990), no. 12, 831--834.

\bibitem{MR} J. C. McConnell and J. C. Robson, {\it Noncommutative Noetherian Rings}, Wiley--Interscience, Chichester, 1987.

\bibitem{NTY} B. Nguyen, K. Trampel and M. Yakimov, Noncommutative discriminants via Poisson primes, arXiv:1603.02585.

\bibitem{RRZ} M. Reyes, D. Rogalski, and J. J. Zhang, Skew Calabi-Yau algebras and homological identities, {\it Adv. Math.}, {\bf 264} (2014), 308--354.

\bibitem{Richard} L. Richard, Sur les endomorphismes des tores quantiques, {\it Comm. Algebra}, {\bf 30} (2002), no.
11, 5283--5306.

\bibitem{RS} L. Richard and A. Solotar, Isomorphisms between quantum generalized Weyl algebras, {\it J.
Algebra Appl.}, {\bf 5} (2006), no. 3, 271--285.

\bibitem{Rigal} L. Rigal, Spectre de l'algb\'{e}bre de Weyl quantique, {\it Beitr. Algebra Geom.}, {\bf 37}
(1996), 119--148.


\bibitem{Tang} X. Tang, Algebra endomorphisms and derivations of some localized down-up algebras, {\it J. Algebra Appl.}, {\bf 14}, 1550034 (2015).


\bibitem{Tsuchimoto} Y. Tsuchimoto, Endomorphisms of Weyl algebras and p--curvature, {\it Osaka J. Math.}, {\bf 42(2)} (2005), 435--452.

\end{thebibliography}
\end{document}